\PassOptionsToPackage{lowtilde}{url}
\documentclass[onefignum,onetabnum]{siamart171218}

\usepackage[final]{microtype}
\usepackage{lmodern}

\usepackage[ruled,vlined,linesnumbered,algo2e]{algorithm2e}
 \usepackage{multirow}
\crefname{algocf}{Algorithm}{Algorithms}

\usepackage{Local_definitions}
\usepackage{subcaption}
\ifpdf
	\hypersetup{
		pdftitle={},
		pdfauthor={T.~Etling, R.~Herzog, E.~Loayza, and G.~Wachsmuth}
	}
\fi

\begin{document}
\maketitle

\begin{abstract}
	We consider shape optimization problems subject to elliptic partial differential equations.
	In the context of the finite element method, the geometry to be optimized is represented by the computational mesh, and the optimization proceeds by repeatedly updating the mesh node positions.
	It is well known that such a procedure eventually may lead to a deterioration of mesh quality, or even an invalidation of the mesh, when interior nodes penetrate neighboring cells.
	We examine this phenomenon, which can be traced back to the ineptness of the discretized objective when considered over the space of mesh node positions.
	As a remedy, we propose a restriction in the admissible mesh deformations, inspired by the Hadamard structure theorem.
	First and second order methods are considered in this setting.
	Numerical results show that mesh degeneracy can be overcome, avoiding the need for remeshing or other strategies.
	\fenics\ code for the proposed methods is available on GitHub.
\end{abstract}



\section{Introduction}
\label{sec:introduction}

Shape optimization is ubiquitous in the design of structures of all kinds, going from drug eluting stents \cite{Zunino2004} until aircraft wings~\cite{SchmidtSchulzIlicGauger2011} or horn-like structures appearing in devices for acoustic or electromagnetic waves~\cite{UdawalpolaBerggren2008}. 
All of these and many other applications involve the solution $u$ of a partial differential equation (PDE), so the general formulation of shape optimization problems considered here is as follows:
\begin{equation}
	\label{eq:general_shape_optimization_problem}
	\min_{\Omega} j(\Omega,u(\Omega))
	.
\end{equation}
Here $u(\Omega)$ is the solution of the underlying PDE defined on the domain $\Omega$, which is to be optimized.
In the following, we will mainly use the reduced objective $J(\Omega) \coloneqq j(\Omega,u(\Omega))$.

Computational approaches to solving PDE-constrained shape optimization problems usually proceed along the following lines.
First, one derives an expression for the \emph{shape derivative} of the objective w.r.t.\ vector fields which describe the perturbation of the current domain $\Omega$.
The perturbations are carried out either in terms of the perturbation of identity, or the velocity method.
We refer the reader to \cite[Chapters~4.3--4.4]{DelfourZolesio2011} for details.
The shape derivative can be stated either as an expression concentrated on the boundary $\partial \Omega$, or as a volume expression.
The first is due to the Hadamard structure theorem (\cite[Theorem~2.27]{SokolowskiZolesio1992}).
For volume expressions, we refer the reader, for instance, to \cite{LaurainSturm2015:1,HiptmairPaganiniSargheini2015:1}.
Second, the shape derivative, which represents a linear functional on the perturbation vector fields, needs to be converted into a vector field $\shapedef$ itself, often referred to as the \emph{shape gradient}.
This can be achieved by evaluating the Riesz representative of the derivative w.r.t.\ an inner product.
The latter is often chosen as the bilinear form associated with the Laplace-Beltrami operator on $\partial \Omega$, or with the linear elasticity (Lam\'e) system on $\Omega$, see e.g.~\cite{SchmidtSchulzIlicGauger2011,SchulzSiebenborn2016,SchmidtSchulz2010:1,SchmidtSchulz2009:2}.
More sophisticated techniques include quasi-Newton or Hessian-based inner products; see~\cite{EpplerHarbrecht2005,NovruziRoche2000,SchulzSiebenbornWelker2015:1,Schulz2014}.
This perturbation field is then used to update the domain $\Omega$ inside a line search method, where the transformed domain
\begin{equation}
	\label{eq:perturbation_of_identity_transformation}
	\Omega_\alpha = \{ x + \alpha \, \shapedef(x): x \in \Omega \} 
\end{equation}
associated with the step size $\alpha$ is obtained from the perturbation of identity approach.

Using the boundary expression of the shape derivative provides an alternative to the domain transformation approach described in the previous paragraph.
In this case, a normal vector field concentrated on the current domain boundary $\partial \Omega$ is obtained, which provides a descent direction.
In the presence of a PDE, modifications of the domain boundary need to be accompanied by a deformation strategy in the interior (see e.g.~\cite{SchmidtIlicSchulzGauger2011}). 
This two-stage process, however, causes discontinuities in the discrete objective function and thus disturbs the performance of optimization algorithms.   
On the other hand, the use of the volume expression has been shown to exhibit superior regularity and better finite element approximation results compared to the boundary expression (see e.g.~\cite{HiptmairPaganiniSargheini2015:1}).

In any case, while the computation of the shape derivative is either based on the continuous or some discrete formulation of problem \eqref{eq:general_shape_optimization_problem}, the computation of the shape gradient and the subsequent updating steps will always be carried out in the discrete setting.
Typically, the shape $\Omega$ is represented by a computational mesh, and the underlying PDE is solved, e.g., by the finite element method.
The perturbation field $V$ is then expressed as a piecewise linear field, i.e., it is represented in terms of a velocity vector attached to each vertex position.
The domain $\Omega$ is subsequently updated according to \eqref{eq:perturbation_of_identity_transformation} inside a line search procedure.

It has been observed in many publications that this straightforward approach has one major drawback: it often leads to a degeneracy of the computational mesh.
This degeneracy manifests itself in different ways, but mostly through degrading cell aspect ratios, or even mesh nodes entering neighboring cells.
\cite{DoganMorinNochettoVerani2007} for instance observe that 
such mesh distortions impair computations and lead to numerical artifacts.
In practice, both phenomena often lead to a breakdown of computational shape optimization procedures.

In \cref{sec:discrete_objective} we shed some light on this process of mesh destruction.
	We attribute it to a discretization artifact, by which the positions of \emph{all} mesh nodes of a computational mesh have an impact on the discrete solution of the PDE present in the problem.
	This presents optimization routines with an opportunity to shift the mesh nodes in such a way that the discrete solution of the PDE exhibits features which allow further descent in the objective, but at the expense of mesh quality and solution accuracy of the PDE.
	Notice that this issue does not arise in the continuous setting, where the redistribution of material points in the interior of the domain has no effect on the PDE solution and thus on the objective.
	In lack of a better name, we refer to the phenomenon described above as \eqq{spurious descent directions} since they are, indeed, leading us further away from the solution of the continuous problem.

Over the past 10~years, a range of various techniques have been proposed to circumvent this major obstacle in computational shape optimization.
A natural choice is to remesh the computational domain; see for instance \cite{WilkeKokGroenwold2005,MorinNochettoPaulettiVerani2012,Sturm2016,DokkenFunkeJohanssonSchmidt2018_preprint,FepponAllaireBordeuCortialDapogny2018_preprint}.
Remeshing can be carried out either in every iteration or whenever some measure of mesh quality falls below a certain threshold.
Drawbacks of remeshing include the high computational cost and the discontinuity introduced into the history of the objective values.

\cite{BaenschMorinNochetto2005,DoganMorinNochettoVerani2007} describe several techniques such as mesh regularization, space adaptivity, angle control in addition to a semi-implicit Euler discretization for the velocity method, with time adaptivity and backtracking line search.
In a follow-up work, \cite{MorinNochettoPaulettiVerani2012} consider a line search method that aims to avoid mesh distortion due to tangential movements of the boundary nodes, combined with a geometrically consistent mesh modification (GCMM) proposed in \cite{BonitoNochettoPauletti2010}. 
\cite{GiacominiPantzTrabelsi2017} address the issue of spurious descent directions, attributed to discretization errors in the underlying PDE model, via a goal-oriented mesh adaptation approach.
Recently, \cite{IglesiasSturmWechsung2017_preprint} proposed to enforce shape gradients from nearly conformal transformations, which are known to preserve angles and ensure a good quality of the mesh along the optimization process.  

Finally, we mention \cite{SchulzSiebenbornWelker2015:1,SchulzSiebenbornWelker2015:2,SchulzSiebenborn2016}, who advocate the linear elasticity model as the inner product to convert shape derivative into a shape gradient.
In particular in \cite{SchulzSiebenbornWelker2015:2} the authors propose to omit the assembly of interior contributions appearing in the discrete volume expression of the shape derivative.
This approach is  related to but conceptionally different from our idea and no analysis is provided there.
A thorough comparison is provided in \cref{subsec:comparison_SSW2016}.

\subsection*{Our Contribution}

In this paper we propose an approach to avoid spurious descent directions in the course of numerical shape optimization procedures, which is different from all of the above and does not require remeshing.
The main idea is based on the observation that---in the continuous setting---shape gradients are perturbation fields which are generated exclusively by normal forces on the boundary of the current domain.
This follows from the Hadamard structure theorem.
However, in the discrete setting, the Hadamard structure theorem is not available, and thus classical discrete shape gradients also contain contributions from interior forces and tangential boundary forces.
We therefore propose to project the shape gradient onto the subspace of perturbation fields generated by normal forces.
We refer to this approach as \emph{restricted mesh deformations}.

We demonstrate that the proposed approach indeed avoids spurious descent directions and degenerate meshes.
As a consequence, we can solve discrete shape optimization problems to high accuracy, i.e., very small norm of the restricted gradient.
Both gradient and Newton schemes in 2D and 3D are considered.
An implementation in the finite element software \fenics\ is available as open-source on GitHub; see \cite{EtlingHerzogLoayzaWachsmuth2018:2}.

The paper is structured as follows.
In \cref{sec:preliminaries} we present a shape optimization model problem and prove, as an auxiliary result, the existence of a globally optimal domain.
In \cref{sec:shape_calculus} we review the volume and boundary representations of the shape derivative.
In \cref{sec:discrete_objective} we consider the discrete counterpart of the model problem and its shape derivative.
We also illustrate the detrimental effect of spurious descent directions.
The main idea of restricted mesh deformations is introduced in \cref{sec:restricted_mesh_deformations}.
An associated \emph{restricted gradient scheme} is also introduced and its performance is compared to the classical shape gradient method in \cref{sec:numerical_results_gradient}.
\Cref{sec:restricted_shape_hessian,sec:numerical_results_newton} are devoted to second-order shape derivatives in the restricted setting and the demonstration of the associated Newton scheme.
Conclusions are given in \cref{sec:conclusions}.

We wish to point out that the model problem considered throughout the paper is clearly academic.
	It was chosen since, as an auxiliary result of the present paper, we present a new technique to prove the existence of optimal shapes, which requires certain properties of the objective to hold; see \cref{rem:more_general_existence}.
	It should be understood that our main idea of considering restricted mesh deformations to avoid spurious descent directions applies to a much broader class of PDE constrained shape optimization problems.

\section{Preliminaries}
\label{sec:preliminaries}

Throughout the paper, we consider the following model problem,
\begin{equation}
	\label{eq:shape_optimization_problem_continuous}
	\text{Minimize} \quad \int_\Omega u \, \dx 
	\quad
	\text{s.t.} \quad  
	\text{$\Omega \subset D$ is open},
	\left\{
		\begin{aligned}
			- \laplace u & = f & & \text{in } \Omega, \\
			u & = 0 & & \text{on } \partial\Omega.
		\end{aligned}
	\right.
\end{equation}
Here the optimization variable $\Omega \subset \R^d$ is an admissible domain contained in some bounded and open hold-all $D \subset \R^d$, and $f \in H^1(D)$ is a given right hand side.
The elliptic state equation is understood in weak form,
\begin{equation}
	\label{eq:weak_form_example_PDE}
	\text{Find } u \in H_0^1(\Omega)
	\qquad\text{such that}\qquad
	\int_\Omega \nabla u \cdot \nabla v \, \dx = \int_\Omega f \, v \, \dx
	\quad\forall v \in H_0^1(\Omega).
\end{equation}

The next result shows that our shape optimization problem
\eqref{eq:shape_optimization_problem_continuous} has a solution
if we slightly relax the class of admissible sets.
We will see that it is sufficient to consider
\emph{quasi-open} rather than open sets.
For an introduction of quasi-open sets, quasi-continuity,
quasi-everywhere (q.e.)
and related notions, we refer the reader to
\cite[Section~5.8]{AttouchButtazzoMichaille2014}.
We consider the slightly relaxed problem
\begin{equation}
	\label{eq:shape_optimization_problem_continuous_relaxed}
	\text{Minimize} \quad \int_\Omega u \, \dx
	\quad
	\text{s.t.} \quad
	\text{$\Omega \subset D$ is quasi-open},
	- \laplace u = f \text{ in } H^{-1}(\Omega)
	.
\end{equation}
Let us recall that
$H_0^1(\Omega) = \{u \in H_0^1(\R^d) \mid u = 0 \text{ q.e.\ in } \R^d\setminus\Omega \}$
and
$H^{-1}(\Omega)$ is the dual space of $H_0^1(\Omega)$.
The PDE in \eqref{eq:shape_optimization_problem_continuous_relaxed} is also to be understood
in the weak sense, i.e.,
\begin{equation*}
	\text{Find } u \in H_0^1(\Omega)
	\qquad\text{such that}\qquad
	\int_D \nabla u \cdot \nabla v = \int_D f \, v \,\dx 
	\quad\forall v \in H_0^1(\Omega).
\end{equation*}
We emphasize that the main reason for this existence result
is that the objective is monotone w.r.t.\ the state $u$,
see also \cref{rem:more_general_existence} below.
\begin{theorem}
	\label{theorem:existence}
	Problem \eqref{eq:shape_optimization_problem_continuous_relaxed}
	admits a global minimizer $(\hat\Omega, \hat u)$.
\end{theorem}
Note that the extreme case $(\hat\Omega, \hat u) = (\emptyset, 0)$
is possible.
\begin{proof}
	First, we remark that it is sufficient to consider
	only pairs $(\{u < 0\}, u)$ with $u \le 0$
	in \eqref{eq:shape_optimization_problem_continuous_relaxed}.
	Indeed,
	if $(\Omega,u)$ is any admissible pair,
	we can consider $(\{u < 0\}, \min(u,0))$ in its stead.
	Note that $\{u < 0\}$ is quasi-open since $u$ can chosen to be quasi-continuous.
	This pair is again admissible due to
	\begin{equation*}
		\int_D \nabla \min(u,0) \cdot \nabla v \, \dx
		=
		\int_\Omega \nabla u \cdot \nabla v \, \dx
		=
		\int_D f \, v \, \dx
		\qquad\forall v \in H_0^1(\{u < 0\}),
	\end{equation*}
	since $v = 0$ q.e.\ on $\Omega \setminus \{u < 0\}$.
	Moreover, the objective value of $(\{u < 0\}, \min(u,0))$
	is not larger than the objective value
	of $(\Omega, u)$.

	Now, let $\{(\Omega_n, u_n)\}$ be a minimizing sequence for \eqref{eq:shape_optimization_problem_continuous_relaxed} with $u_n \le 0$
	and $\Omega_n = \{u_n < 0\}$.
	It is clear that the sequence $\{u_n\}$
	is bounded in $H_0^1(D)$,
	therefore we can extract a weakly convergent subsequence
	(without relabeling)
	with weak limit $u$.
	Clearly, $u \le 0$.
	Now we define $\hat\Omega = \{u < 0\}$
	and denote by $\hat u \in H_0^1(\hat\Omega)$
	the solution of $-\Delta \hat u = f$ in $H^{-1}(\hat\Omega)$.
	It remains to check that $\hat u \le u$ holds
	since this implies the global optimality of $\hat u$
	(due to the monotonicity of the objective).
	To this end, we
	choose an arbitrary $v \in H_0^1(D)$ such that
	$-u \ge v \ge 0$.
	For $v_n := \min(-u_n, v)$
	we have $v_n \in H_0^1(\Omega_n)$
	due to $v \ge 0$.
	Moreover, $v_n \weakly \min(-u,v) = v$ in $H_0^1(D)$, see \cite[Lemma~4.1]{Wachsmuth2014:2}.
	Thus,
	\begin{align*}
		\int_D f \, v \, \dx
		&=
		\lim_{n\to\infty} \int_D f \, v_n \, \dx
		=
		\lim_{n\to\infty} \int_D \nabla u_n \cdot \nabla v_n \, \dx
		\\&
		=
		\lim_{n\to\infty} \int_D
		\nabla (u_n + v) \cdot \nabla(v_n - v)
		+ \nabla u_n \cdot \nabla v
		- \nabla v \cdot \nabla (v_n - v)
		\, \dx
		\\&
		=
		\lim_{n\to\infty} \int_D -\nabla\abs{\min(-u_n-v,0)}^2 + \nabla u_n \cdot \nabla v \, \dx
		\le
		\int_D \nabla u \cdot \nabla v \, \dx
		.
	\end{align*}
	Since $v \in H_0^1(\hat\Omega)$, we can test the equation for $\hat u$ with $v$ and we find
	\begin{equation*}
		\int_{\hat\Omega} \nabla(\hat u - u) \cdot \nabla v  \, \dx \le 0
		\qquad\forall v \in H_0^1(D) \text{ satisfying } {-u} \ge v \ge 0.
	\end{equation*}
	Now, we can use a density argument,
	see \cite[Lemme~3.4]{Mignot1976},
	to obtain that this inequality holds for all
	$v \in H_0^1(\hat\Omega)$ which satisfy $v \ge 0$.
	Using $v = \max(\hat u - u, 0)$ implies $\max(\hat u - u, 0) = 0$,
	i.e., $\hat u \le u$.
	Finally, the optimality of $(\hat\Omega, \hat u)$ follows from
	\begin{equation*}
		\int_D \hat u \, \dx
		\le
		\int_D u \, \dx
		=
		\lim_{n \to \infty} \int_D u_n \, \dx.
	\end{equation*}
\end{proof}
\begin{remark}
	\label{rem:capacitary_measures}
	There is a deeper reason for $\hat u \le u$ being true in the above proof.
	Indeed, using the theory of relaxed Dirichlet problems,
	one can show that $u$ satisfies $-\Delta u + \mu \, u = f$
	for some capacitary measure $\mu$.
	We refer to \cite[Section~5.8.4]{AttouchButtazzoMichaille2014}
	for a nice introduction to capacitary measures.
	Due to $u \le 0$ we have
	(in a certain sense) $\mu \, u \le 0$ and therefore
	$\hat u \le u$ follows from the maximum principle since
	``$-\Delta \hat u = f \le f - \mu\,u = -\Delta u$''.
	However, we included the above direct proof because it does not rely on the notion
	of capacitary measures.
\end{remark}
\begin{remark}
	\label{rem:more_general_existence}
	The above proof of existence generalizes to a larger class of objective functionals.
	In fact, we can replace the objective in \eqref{eq:shape_optimization_problem_continuous_relaxed}
	with
	\begin{equation*}
		\int_\Omega j(x, u(x)) \, \dx
	\end{equation*}
	if the integrand $j$ satisfies
	\begin{subequations}
		\label{eq:assumptions_integrand}
		\begin{align}
			& j(x, \cdot) \text{ is monotonically increasing on $(-\infty,0]$ and non-negative on $[0,\infty)$}
			,\\
			& j(\cdot, u) \in L^1(D) \quad \forall u \in H_0^1(D)
			,\\
			& u_n \weakly u \text{ in $H_0^1(D)$ implies } \int_D j(u) \, \dx \le \liminf_{n\to\infty} \int_D j(u_n) \, \dx
			.
		\end{align}
	\end{subequations}
	Under these general assumptions, one can use the same proof as the one given for \cref{theorem:existence} above,
	but the final estimate has to be replaced by
	\begin{align*}
		\int_{\hat\Omega} j(\cdot, \hat u) \, \dx
		\le
		\int_{\hat\Omega} j(\cdot, u) \, \dx
		&=
		\int_D j(\cdot, u) - j(\cdot, 0) \, \dx + \int_{\{u < 0\}} j(0) \, \dx
		\\
		&\le
		\liminf_{n \to \infty}
		\int_D j(\cdot, u_n) - j(\cdot, 0) \, \dx + \int_{\{u_n < 0\}} j(0) \, \dx
		\\&
		=
		\liminf_{n \to \infty}
		\int_{\Omega_n} j(\cdot, u_n) \, \dx
		.
	\end{align*}
	Note that Fatou's lemma together with $u_n \to u$ a.e.\ (along a subsequence)
	implies
	\begin{equation*}
		\int_{\{u < 0\}} j(0) \, \dx
		\le
		\liminf_{n \to \infty}
		\int_{\{u_n < 0\}} j(0) \, \dx
		.
	\end{equation*}
	Again, this shows the optimality of $(\hat\Omega, \hat u)$.
\end{remark}

\section{Shape Calculus}
\label{sec:shape_calculus}

This section is devoted to the presentation of the shape differentiability of problem \eqref{eq:shape_optimization_problem_continuous}.
Since this is rather standard problem we will be able to directly apply results from \cite{ItoKunischPeichl2008:1}.
To this end, we assume that both the hold-all $D \subset \R^d$ and $\Omega \subset \R^d$ are open and have $C^{1,1}$-boundaries $\partial D$ and $\partial \Omega$, respectively.
Moreover we assume $\overline{\Omega} \subset D$ so that $\Omega$ has a positive distance to the boundary of $D$.

We are describing variations of the domain $\Omega$ by the \emph{perturbation of identity method}, i.e., we consider a family of transformations $\{T_\alpha\}_{\alpha\in [0, \tau]}$ such that 
\begin{equation}
	\label{eq:perturbation_by_id}
	T_\alpha = \id + \alpha \, V,
\end{equation} 
where $V \in C^{1,1}(D)^d$ is a given vector field.
The family $\{T_\alpha\}$ creates a family of perturbed domains $\Omega_\alpha = T_\alpha(\Omega)$.
In view of Banach's fixed point theorem, there exists a bound $\tau > 0$ such that $T_\alpha$ is invertible for all $\alpha \in [0,\tau]$.

By a straightforward application of \cite[Theorem~2.1]{ItoKunischPeichl2008:1} we obtain the following result.
\begin{theorem}
	\label{theorem:existence_shape_derivative}
	The shape functional given in~\eqref{eq:shape_optimization_problem_continuous} is shape differentiable and its shape derivative in the direction of the perturbation field $V$ is given by 
	\begin{multline}
		\label{eq:shape_derivative}
		J'(\Omega;V) 
		= 
		\int_{\Omega} u \, (\div V) \, \dx 
		\\
		+ 
		\int_\Omega (\nabla u)^\top \, \bigh[]{ (\div V) \, \id - DV - DV^\top } \, \nabla p \, \dx
		- 
		\int_\Omega \div (f \, V) \, p \, \dx 
	\end{multline}
	where $DV$ denotes the Jacobian of $V$ and the adjoint state $p$ is the unique solution of the following adjoint problem,
	\begin{equation}
		\label{eq:adjoint_problem}
		\text{Find } p \in H_0^1(\Omega)
		\quad \text{such that} \quad
		\int_D \nabla p \cdot \nabla v \, \dx = - \int_D  v \, \dx 
		\quad \text{for all } v \in H_0^1(\Omega).
	\end{equation}
\end{theorem}

Notice that \eqref{eq:shape_derivative} is the so-called volume or weak formulation of the shape derivative of \eqref{eq:shape_optimization_problem_continuous}.
Besides the volume formulation, there exists an alternative representation of \eqref{eq:shape_derivative} by virtue of the well known Hadamard structure theorem; see \cite[Chapter~9, Theorem~3.6]{DelfourZolesio2011}.
We state it here in a particularized version for problem \eqref{eq:shape_optimization_problem_continuous}. 
From now on, $\nu$ denotes the outer unit normal vector along the boundary $\partial \Omega$ of $\Omega$.

\begin{corollary}[Hadamard structure theorem for \eqref{eq:shape_optimization_problem_continuous}]
	\label{corollary:Hadamard}
	The shape derivative \eqref{eq:shape_derivative} of problem \eqref{eq:shape_optimization_problem_continuous} has the representation 
	\begin{equation}
		\label{eq:shape_derivative_Hadamard}
		J'(\Omega;V) = 
		\int_{\partial \Omega} g_\Omega \, (V \cdot \nu) \, \ds
		\quad
		\text{with } 
		g_\Omega = - \frac{\partial u}{\partial \nu}\,\frac{\partial p}{\partial \nu}.
	\end{equation}
\end{corollary}
Notice that under the assumption that $\Omega$ has a $C^{1,1}$-boundary, $u$ and $p$ belong to $H^2(\Omega)$ and thus their normal derivatives are in $H^{1/2}(\partial \Omega)$, which embeds into $L^4(\partial \Omega)$ when $d \le 3$; see for instance \cite[Theorem~4.12]{AdamsFournier2003}.
Consequently, $g_\Omega = - \frac{\partial u}{\partial \nu}\,\frac{\partial p}{\partial \nu}$ belongs to $L^2(\partial \Omega)$ in this case.

Formula \eqref{eq:shape_derivative_Hadamard} is known as the boundary or strong representation of \eqref{eq:shape_derivative}, and it can be obtained from \eqref{eq:shape_derivative} by the divergence theorem; compare~\cite{Sturm2015:1}, \cite[Chapter~3.3]{SokolowskiZolesio1992}, \cite[Example~3.3]{HaslingerMaekinen2003}. 
We also refer the reader to \cite{HiptmairPaganiniSargheini2015:1}, where the volume and boundary formulations are compared w.r.t.\ their order of convergence in a finite element setting.

\section{Investigation of the Discrete Objective}
\label{sec:discrete_objective}

In order to solve the shape optimization problem \eqref{eq:shape_optimization_problem_continuous} numerically, some kind of discretization has to be applied.
The most common choice in the literature consists in a discretization of the PDE by some finite element space defined over a computational mesh, which we denote by $\Omega_h$ and whose nodal positions serve to represent the discrete unknown domain.

A common choice is to replace $H^1_0(\Omega)$ by the finite element space of piecewise linear, globally continuous functions,
\begin{equation}
	\label{eq:FE_space_example_PDE}
	S_0^1(\Omega_h) = \{ u \in H^1_0(\Omega_h): \restr{u}{T} \in \PP_1(T) \text{ for all cells $T$ in $\Omega_h$} \}
\end{equation}
defined over an approximation $\Omega_h$ of $\Omega$ consisting of geometrically conforming simplicial cells, i.e., triangles and tetrahedra in $d = 2$ or $d = 3$ space dimensions, respectively.
Consequently, the state equation \eqref{eq:weak_form_example_PDE} is replaced by
\begin{equation}
	\label{eq:weak_form_example_PDE_discrete}
	\text{Find } u_h \in S_0^1(\Omega_h)
	\quad\text{such that}\quad
	\int_\Omega \nabla u_h \cdot \nabla v_h \, \dx = \int_\Omega f \, v_h \, \dx
	\quad \forall v_h \in S_0^1(\Omega_h).
\end{equation}
This leads to the following discrete version of \eqref{eq:shape_optimization_problem_continuous} frequently encountered in the literature,
\begin{equation}
	\begin{aligned}
		\label{eq:shape_optimization_problem_discrete}
		\text{Minimize} \quad & \int_{\Omega_h} u_h \, \dx \quad \text{w.r.t.\ } u_h \in S_0^1(\Omega_h) \text{ and the nodal positions in } \Omega_h
		\\
		\text{s.t.} \quad & \eqref{eq:weak_form_example_PDE_discrete}
		.
	\end{aligned}
\end{equation}
We refer the reader to~\cite{GanglLangerLaurainMeftahiSturm2015:1_preprint,Sturm2016,SchulzSiebenbornWelker2015:2, SchulzSiebenborn2016} for examples of this procedure.

Let us denote by $J_h(\Omega_h)$ the reduced objective value in \eqref{eq:shape_optimization_problem_discrete}, i.e., $J_h(\Omega_h) = \int_{\Omega_h} u_h \, \dx$, where $u_h$ is the unique solution of \eqref{eq:weak_form_example_PDE_discrete}.
In order to derive a discrete variant of the volume formulation \eqref{eq:shape_derivative} of the shape derivative, we introduce the discrete adjoint equation,
\begin{equation}
	\label{eq:adjoint_problem_discrete}
	\text{Find } p_h \in S_0^1(\Omega_h)
	\quad \text{such that} \quad
	\int_{\Omega_h} \nabla p_h \cdot \nabla v_h \, \dx = - \int_{\Omega_h}  v_h \, \dx 
	\quad \text{for all } v_h \in S_0^1(\Omega_h).
\end{equation}

The following theorem shows that a straightforward replacement of the state $u$ and adjoint state $p$ by their finite element equivalents $u_h$ and $p_h$ in \eqref{eq:shape_derivative} yields the correct formula for the shape derivative $J_h'(\Omega_h;\shapedefdisc)$ of the discrete objective $J_h$, provided that the perturbation field $\shapedefdisc$ is piecewise linear, i.e., $\shapedefdisc$ belongs to
\begin{equation}
	\label{eq:FE_space_deformations}
	S^1(\Omega_h)^d = \{ u \in H^1(\Omega_h)^d: \restr{u}{T} \in \PP_1(T)^d \text{ for all cells $T$ in $\Omega_h$} \}
	.
\end{equation}

\begin{theorem}
	\label{theorem:discrete_shape_derivative}
	Suppose that $u_h$ and $p_h$ are the unique weak solutions of the discrete state equation \eqref{eq:weak_form_example_PDE_discrete}, and the discrete adjoint equation \eqref{eq:adjoint_problem_discrete}, respectively.
	Moreover, let $\shapedefdisc \in S^1(\Omega_h)^d$.
	Then 
	\begin{multline}
		\label{eq:shape_derivative_discrete}
		J_h'(\Omega_h;\shapedefdisc) 
		= 
		\int_{\Omega_h} u_h \, (\div \shapedefdisc) \, \dx 
		\\
		+ 
		\int_{\Omega_h} (\nabla u_h)^\top \, \bigh[]{ (\div \shapedefdisc) \, \id - D\shapedefdisc - D\shapedefdisc^\top } \, \nabla p_h \, \dx 
		- 
		\int_{\Omega_h} \div (f \, \shapedefdisc) \, p_h \, \dx 
		.
	\end{multline}
\end{theorem}
The proof of this theorem follows along the lines of the continuous case, see, e.g., \cite{HiptmairPaganiniSargheini2015:1,LaurainSturm2015:1}.
A detailed derivation can be found in \cite[Section~4]{DelfourPayreZolesio1985}.

\begin{remark}
	\label{remark:shape_derivative_discrete}
	\begin{enumerate}
		\item 
			\cref{theorem:discrete_shape_derivative} can be viewed as the statement that discretization and optimization (in the sense of forming the shape derivative) commute for problem \eqref{eq:shape_optimization_problem_continuous}.

		\item 
			The finite element analogue of the \emph{boundary expression} \eqref{eq:shape_derivative_Hadamard} is \emph{not} an exact representation of the discrete shape derivative.
			This is since the integration by parts necessary to pass from the volume to the boundary expression has to be done element by element and it leaves inter-element contributions; see also the discussion in \cite{Berggren2010}.

		\item 
			\cref{theorem:discrete_shape_derivative} remains true when higher order Lagrangian finite elements on simplices are used in place of $S_0^1(\Omega_h)$.
			However it is essential that $\shapedefdisc$ remains piecewise linear so piecewise polynomials are transformed into piecewise polynomials of the same order.

		\item 
			Alternative expressions for \eqref{eq:shape_derivative_discrete} can be obtained following the so-called discrete adjoint approach, in which the derivative of $J_h(\Omega_h)$ w.r.t.\ the nodal positions of $\Omega_h$ is addressed by differentiating the finite element matrices.
			We refer to \cite{SchneiderJimack2008,Berggren2010,RothUlbrich2013} for examples of this procedure.
	\end{enumerate}
\end{remark}

Despite the simplicity to obtain the shape derivative of the discrete problem, we would like to emphasize here that the discrete problem \eqref{eq:shape_optimization_problem_discrete} itself has the following serious drawback.
The search space obtained from utilizing the nodal positions of the mesh $\Omega_h$ as optimization variables includes meshes with very degenerate cells.
Those lead to poor approximations of solutions of the state equation, which may give rise, however, to smaller values of the discrete objective.
Therefore, any optimization algorithm for the solution of \eqref{eq:shape_optimization_problem_discrete} sooner or later is likely to encounter spurious descent directions which typically have support in only a few mesh nodes and which lead to degenerate meshes.

\begin{figure}
	\includegraphics[width=0.33\linewidth]{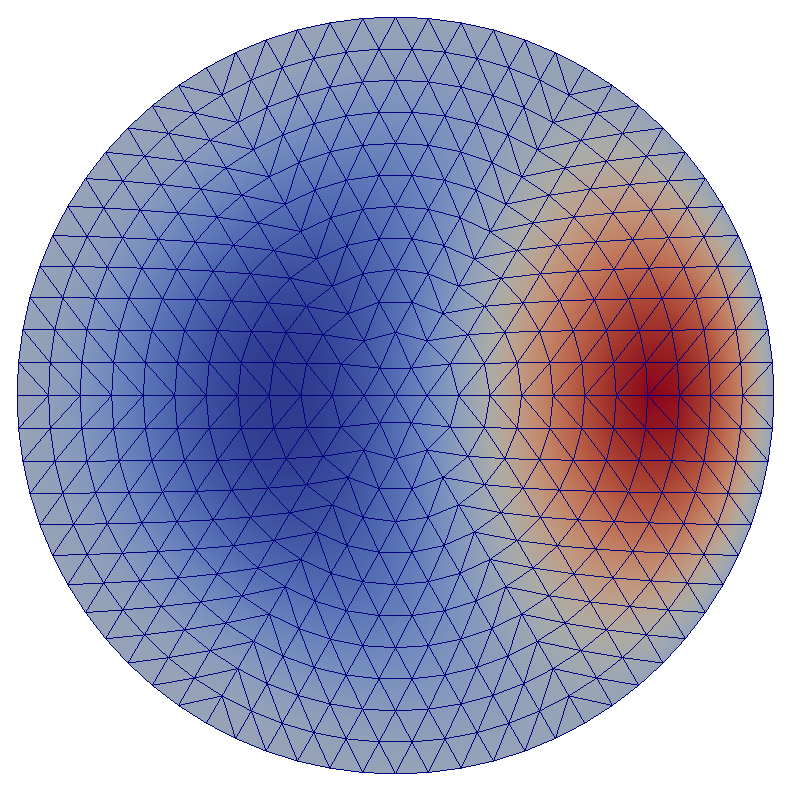}\hfill
	\includegraphics[width=0.33\linewidth]{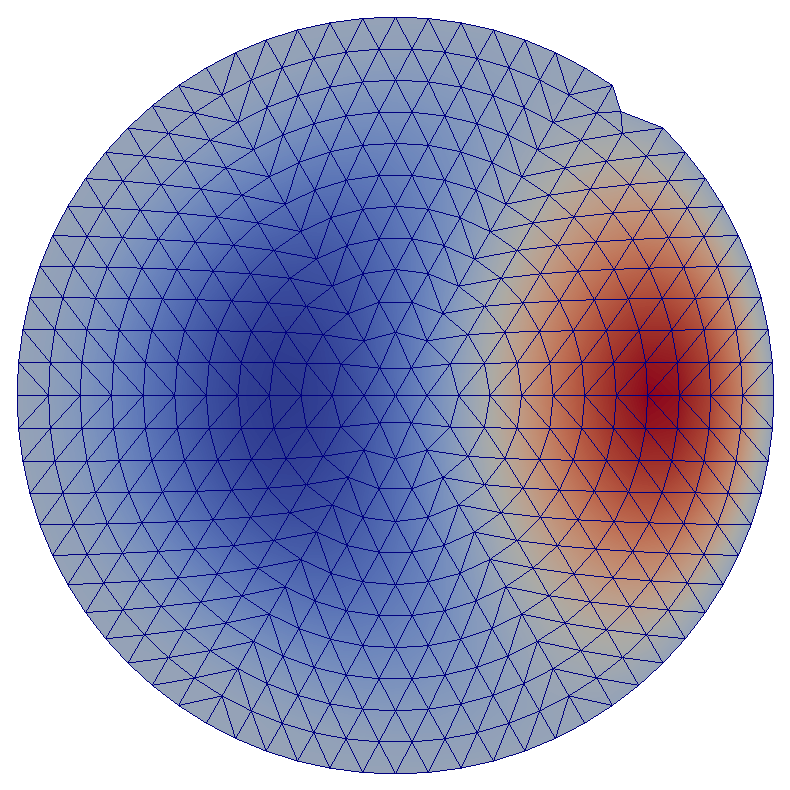}\hfill
	\includegraphics[width=0.33\linewidth]{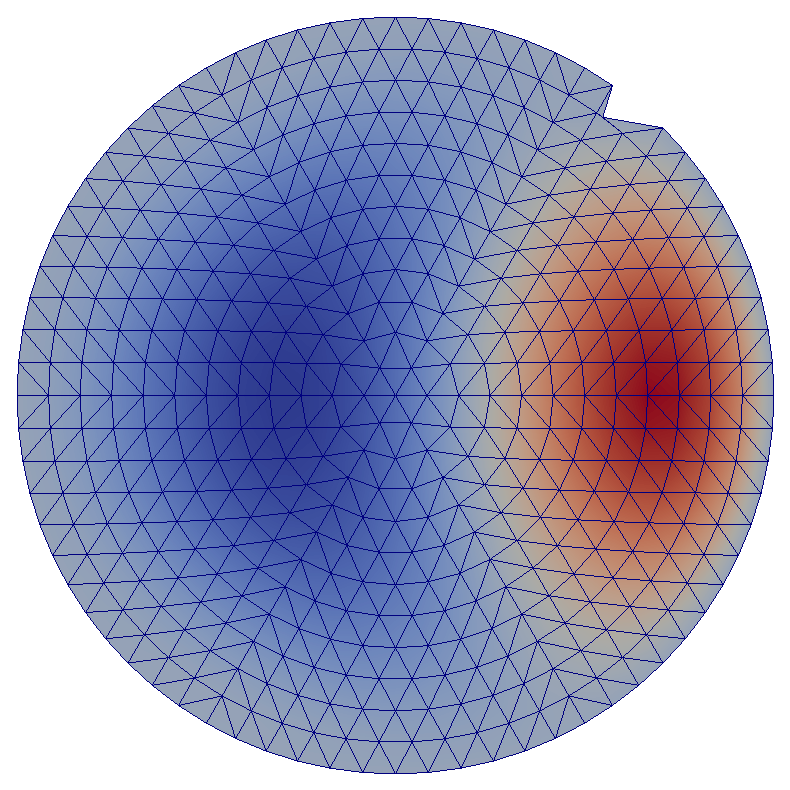}
	\caption{Evolution of the mesh under $\alpha \mapsto T_\alpha (\Omega_h)$ with perturbation field $\shapedefdisc$ given in \cref{example:spurious_descent_direction} at $\alpha = 0.00$, $\alpha = 0.05$, $\alpha = 0.10$ (from left to right).
	The solutions $u_h$ of the state equation \eqref{eq:weak_form_example_PDE_discrete} are also shown.}
	\label{fig:destroy_mesh}
\end{figure}

\begin{example}
	\label{example:spurious_descent_direction}
	Let us illustrate this behavior by means of problem \eqref{eq:shape_optimization_problem_continuous} with data $f(x,y) = 2.5 \, (x + 0.4 - y^2)^2 + x^2 + y^2 - 1 $.
	The optimal domain $\Omega$ is unknown.
	We begin with the computational mesh $\Omega_h$ shown in \cref{fig:destroy_mesh}~(left).
	Consider for example the piecewise linear vector field $\shapedefdisc$ represented by its nodal values
	\begin{equation*}
		\shapedefdisc = 
		\begin{cases}
			(-0.9510,-0.3090)^\top, & \text{for the node } v_0,\\
			(0,0)^\top, & \text{for all other nodes}
		\end{cases}
	\end{equation*}
	where the boundary node $v_0$ can be easily identified from \cref{fig:destroy_mesh}.
\end{example}

\begin{figure}
	\centering
	\includegraphics[scale=1]{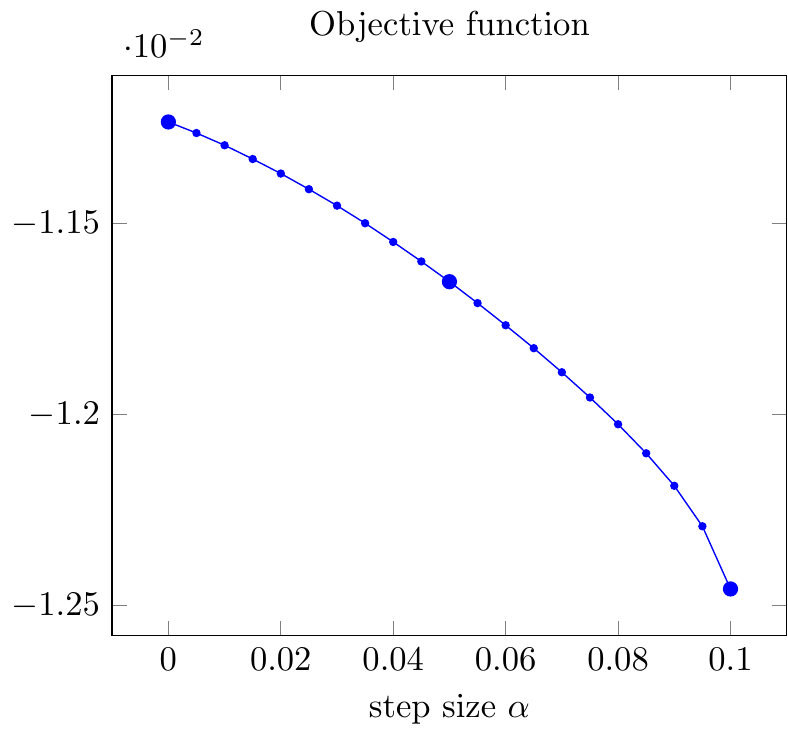}
	\caption{$\alpha \mapsto J\bigh(){ T_\alpha (\Omega_h) }$ for \cref{example:spurious_descent_direction}. The step sizes $\alpha = 0.00$, $\alpha = 0.05$, $\alpha = 0.10$ belonging to the domains in \cref{fig:destroy_mesh} are highlighted.}
	\label{fig:destroy_mesh_decay}
\end{figure}

We found that $V$ is not only a descent direction for the objective at $\Omega_h$ but in fact that the line search function
\begin{equation*}
	\alpha \mapsto J\bigh(){ T_\alpha (\Omega_h) },
	\quad
	T_\alpha = \id + \alpha \, \shapedefdisc
\end{equation*}
decreases until the triangle formed by $v_0$ and its two interior neighbors degenerates to a line, which happens at $\alpha = 0.1$; see \cref{fig:destroy_mesh_decay}.
At this point, finite element computations break down.

In computational experience spurious descent directions do not usually occur during the early iterates.
Thus they can be, and often are, avoided by early stopping, at the expense of a reduced tolerance.
Alternatively, mesh quality control and remeshing can help to avoid mesh destruction, but this introduces discontinuities in the objective function's history.

In any case, the existence of spurious descent directions is a structural disadvantage of problem \eqref{eq:shape_optimization_problem_discrete}.
Therefore we propose in the following section a new computational approach.
Our approach does not seek to solve \eqref{eq:shape_optimization_problem_discrete} literally but in a certain relaxed sense, which is inspired by the Hadamard structure theorem and which avoids spurious descent directions.

\section{Restricted Mesh Deformations}
\label{sec:restricted_mesh_deformations}

By the Hadamard structure theorem, 
the shape derivative 
for the continuous problem
consists of normal boundary forces only,
see \eqref{eq:shape_derivative_Hadamard} above.
This is no longer the case for the discrete problem.
The reason is that the finite element solutions $u_h$ and $p_h$ are only of limited regularity,
and thus a global integration by parts necessary to pass from the volume expression \eqref{eq:shape_derivative_discrete} to a boundary expression is not available.
This has been pointed out, for instance, in \cite[note on p.~562]{DelfourZolesio2011}.
Therefore, we are going to continue with the discretely exact volume expression \eqref{eq:shape_derivative_discrete} but mimic the behavior of the continuous setting in the evaluation of the shape gradient, where we alloy only for shape displacements which are induced by normal forces.

\subsection{Continuous Setting}
\label{subsec:restricted_mesh_deformations_continuous}
To illustrate the situation,
we start by discussing the continuous case.
We have seen in
\eqref{eq:shape_derivative}
that the shape derivative $J'(\Omega;\cdot)$
is an element of a dual space, e.g.\ an element of $(W^{1,\infty}(\Omega)^d)\dualspace$.
In order to utilize this information
for moving the domain $\Omega$,
we have to convert this dual element
into a proper function.
We follow the approach of
\cite{SchulzSiebenbornWelker2015:2}.
To this end, we introduce the elasticity operator
$E : H^1(\Omega)^d \to (H^1(\Omega)^d)\dualspace$
via
\begin{equation}
	\label{eq:elasticity_operator}
	\dual{E\,\shapedef}{W}
	:=
	\int_\Omega
	2 \, \mu \, \bvarepsilon(\shapedef) \dprod \bvarepsilon(W)
	+
	\lambda \, \trace(\bvarepsilon(\shapedef)) \, \trace(\bvarepsilon(W))
	+
	\delta \, \shapedef \cdot W
	\,\dx
\end{equation}
for all $\shapedef, W \in H^1(\Omega)^d$.
Here and throughout, $D$ denotes the derivative (Jacobian) of a vector valued function, $\bvarepsilon(\shapedef) = (D\shapedef + D\shapedef^\top)/2$ is the linearized strain tensor,
$\mu, \lambda$ are the Lamé parameters
and $\delta > 0$ is a damping term.
We assume $\mu > 0$, $d \, \lambda + 2\,\mu > 0$ so that $E$ becomes positive semi-definite on $H^1(\Omega)^d$.
Note that
we do not consider Dirichlet boundary conditions in the space $H^1(\Omega)^d$.
Therefore a positive damping parameter $\delta > 0$ is needed to ensure the coercivity of $E$,
i.e., $\dual{E\,\shapedef}{\shapedef} \ge \underline c \, \norm{\shapedef}_{H^1(\Omega)^d}^2$ 
with some $\underline c > 0$.
This result is due to Korn's inequality, see for instance \cite[Proposition~6.6.1]{AttouchButtazzoMichaille2014}.
Thus, $E$ is an isomorphism and it furnishes $H^1(\Omega)^d$ with an inner product $\scalarprod{\shapedef}{W}_E \coloneqq \dual{E\,\shapedef}{W}$ so that $E$ becomes the associated Riesz isomorphism.

In order to avoid technical regularity issues,
we assume that the shape derivative
\eqref{eq:shape_derivative}
enjoys the higher regularity
$J'(\Omega; \cdot) \in (H^1(\Omega)^d)\dualspace$.
This holds, e.g., if $\Omega$ is sufficiently smooth,
due to the higher regularity of $u$ and $p$.
In order to compute the negative shape gradient w.r.t.\ the $E$-inner product
on the continuous level,
we solve
\begin{equation}
	\label{eq:continuous_shape_gradient_1}
	\text{Minimize}\quad
	J'(\Omega; \shapedef)
	+
	\frac12 \, \dual{E \shapedef}{\shapedef}
	\quad \text{s.t.\ } \shapedef \in H^1(\Omega)^d.
\end{equation}
The solution of this problem yields the negative shape gradient
\begin{equation}
	\label{eq:shape_gradient}
	\shapegrad := -E^{-1} J'(\Omega; \cdot).
\end{equation}
Now,
we introduce the normal force operator
$N : L^2(\partial\Omega) \to (H^1(\Omega)^d)\dualspace$
given by
\begin{equation}
	\label{eq:normal_force_operator}
	\dual{N F}{\shapedef}
	=
	\int_{\partial\Omega} F \, (\shapedef \cdot \nu) \, \ds
\end{equation}
for all
$F \in L^2(\partial\Omega)$ and $\shapedef \in H^1(\Omega)^d$.
Using again \eqref{eq:shape_derivative_Hadamard},
we find that $J'(\Omega; \cdot)$ can be written as
$J'(\Omega; \cdot) = N g_\Omega$ with
\begin{equation*}
	g_\Omega
	=
	- \frac{\partial u}{\partial \nu}\,\frac{\partial p}{\partial \nu}
	\in L^2(\partial\Omega).
\end{equation*}
Therefore, it is easy to see that
problem \eqref{eq:continuous_shape_gradient_1}
is equivalent to
\begin{equation}
	\label{eq:continuous_shape_gradient_2}
	\begin{aligned}
		\text{Minimize}\quad &
		J'(\Omega; \shapedef)
		+
		\frac12 \, \dual{E \shapedef}{\shapedef}
		\\
		\text{with respect to}\quad&
		\shapedef \in H^1(\Omega)^d,
		F \in L^2(\partial\Omega)
		\\
		\text{such that}\quad&
		E \, \shapedef - N \, F = 0
		.
	\end{aligned}
\end{equation}
Indeed, the additional constraint $E \, \shapedef - N \, F = 0$ is automatically satisfied by the unconstrained solution of \eqref{eq:continuous_shape_gradient_1}.
However, we will see
that this property is lost
after discretization,
i.e., the discrete counterparts
of \eqref{eq:continuous_shape_gradient_1} and \eqref{eq:continuous_shape_gradient_2}
are going to differ.
Note that the solution $(V,F)$ of \eqref{eq:continuous_shape_gradient_2} is unique
due to coercivity of $E$ and injectivity of $N$.
Moreover, since $\begin{bmatrix} E & -N \end{bmatrix}$ is surjective, there exists a unique Lagrange multiplier $\Pi \in H^1(\Omega)^d$ associated with the constraint $E\,\shapedef - N\,F = 0$; see for instance \cite[Chapter~9.3, Theorem~1]{Luenberger1969}.
We therefore obtain the following necessary and sufficient optimality conditions for \eqref{eq:continuous_shape_gradient_2} 
in saddle-point form,
\begin{equation}
	\label{eq:continuous_saddle_point_system}
	\begin{pmatrix}
		E & 0 & E \\
		0 & 0 & -N\adjoint\\
		E & -N & 0
	\end{pmatrix}
	\begin{pmatrix}
		\shapedef \\ F \\ \Pi
	\end{pmatrix}
	=
	\begin{pmatrix}
		-J'(\Omega; \cdot) \\ 0 \\ 0
	\end{pmatrix}
	.
\end{equation}
Here,
$N\adjoint : H^1(\Omega)^d \to L^2(\partial\Omega)$ is the adjoint of $N$,
where we identified $L^2(\partial\Omega)$ with its dual space.
The multiplier $\Pi$ in \eqref{eq:continuous_saddle_point_system} necessarily satisfies $\Pi = 0$ since $E$ is bijective.
Now,
it is easy to see that \eqref{eq:continuous_saddle_point_system} is equivalent to
solving
\begin{subequations}
	\label{eq:continuous_saddle_point_system_reduced}
	\begin{align}
		\begin{pmatrix}
			0 & N\adjoint\\
			N & E
		\end{pmatrix}
		\begin{pmatrix}
			F \\ \Pi
		\end{pmatrix}
		&=
		\begin{pmatrix}
			0 \\ -J'(\Omega; \cdot)
		\end{pmatrix}
		,
		\\
		\shapedef
		&=
		-E^{-1} J'(\Omega;\cdot) - \Pi
		.
	\end{align}
\end{subequations}
Recall that $-E^{-1} J'(\Omega;\cdot)$
is the usual negative shape gradient w.r.t.\ $E$
(i.e., the solution of \eqref{eq:continuous_shape_gradient_1}),
whereas
$-\Pi$
is a correction
in order to obtain a shape displacement
in the subspace $\im(E^{-1} N)$.
Again,
we emphasize that
we have
$\Pi = 0$
in the continuous setting,
due to $-J'(\Omega;\cdot) = -N \, g_\Omega$.
Therefore,
the solution of \eqref{eq:continuous_saddle_point_system_reduced}
is just
the usual shape gradient
$\shapegrad = -E^{-1} J'(\Omega; \cdot)$.

Before discussing the discretized setting,
we note that
\eqref{eq:continuous_shape_gradient_2}
is equivalent to
\begin{equation}
	\label{eq:continuous_shape_gradient_3}
	\begin{aligned}
		\text{Minimize}\quad &
		\frac12 \, \bigdual{E (\shapedef - \shapegrad)}{\shapedef - \shapegrad}
		\\
		\text{with respect to}\quad&
		\shapedef \in H^1(\Omega)^d,
		F \in L^2(\partial\Omega)
		\\
		\text{such that}\quad&
		E \, \shapedef - N \, F = 0
		.
	\end{aligned}
\end{equation}
Hence, the solution $\shapedef$
is the orthogonal projection (w.r.t.\ the inner product induced by $E$)
of the usual shape gradient $\shapegrad = -E^{-1} J'(\Omega; \cdot)$
into the space $\im(E^{-1} N)$,
i.e., the space of deformations induced by normal forces.
This motivates to denote the solution of \eqref{eq:continuous_shape_gradient_2}
by $\shapeprojgrad$.

\subsection{Discretized Setting}
\label{subsec:restricted_mesh_deformations_discrete}

Next, we discuss the discretized setting.
We refer to \cref{sec:discrete_objective} above
for the introduction of the finite-element discretization.
In addition to the FE space $S_0^1(\Omega_h) \subset H_0^1(\Omega_h)$,
we recall from \eqref{eq:FE_space_deformations} the discrete space of mesh deformations
\begin{equation*}
	S^1(\Omega_h)^d = \{ u \in H^1(\Omega_h)^d: \restr{u}{T} \in \PP_1(T)^d \text{ for all cells $T$ in $\Omega_h$} \}
\end{equation*}
and the boundary space
\begin{equation}
	\label{eq:FE_space_boundary}
	S^1(\partial\Omega_h) = \{ u \in C(\partial\Omega_h)^d: \restr{u}{E} \in \PP_1(E)^d \text{ for all edges $E$ on $\partial\Omega_h$} \}
	.
\end{equation}
We recall that the discrete shape derivative
$J_h'(\Omega_h;\cdot) \in (S^1(\Omega_h)^d)\dualspace$
was given in \eqref{eq:shape_derivative_discrete}.
Moreover, the discretization directly leads to the
discretized operators $E_h : S^1(\Omega_h)^d \to (S^1(\Omega_h)^d)\dualspace$,
$N_h : S^1(\partial\Omega_h) \to (S^1(\Omega_h)^d)\dualspace$
which are defined via
\begin{align*}
	\dual{E_h \, \shapedefdisc}{W_h}
	&:=
	\int_{\Omega_h}
	2 \, \mu \, \bvarepsilon(\shapedefdisc) \dprod \bvarepsilon(W_h)
	+
	\lambda \, \trace(\bvarepsilon(\shapedefdisc)) \, \trace(\bvarepsilon(W_h))
	+
	\delta \, \shapedefdisc \cdot W_h
	\,\dx
	,
	\\
	\dual{N_h \, F_h}{\shapedefdisc}
	&:=
	\int_{\partial\Omega_h} F_h \, (\shapedefdisc \cdot \nu) \, \ds
\end{align*}
for all $\shapedefdisc,W_h \in S^1(\Omega_h)^d$ and $F_h \in S^1(\partial\Omega_h)$.
Next, we will investigate the discrete counterparts of \eqref{eq:continuous_shape_gradient_1}
and \eqref{eq:continuous_shape_gradient_2}.
The straightforward discretization of \eqref{eq:continuous_shape_gradient_1}
reads
\begin{equation}
	\label{eq:discrete_shape_gradient_1}
	\text{Minimize}\quad
	J_h'(\Omega_h; \shapedefdisc)
	+
	\frac12 \, \dual{E_h \shapedefdisc}{\shapedefdisc}.
\end{equation}
We denote its unique solution by $\shapegraddisc$.

The important difference to the continuous case is
that Hadamard's structure theorem is not available.
The reason is that the discrete state $u_h$
has only the limited regularity $u_h \in H^1_0(\Omega_h)$
and this regularity is not enough to
transform the domain integral into a boundary integral
via integration by parts,
see the last paragraph in chapter~10, section~5.6 
of \cite{DelfourZolesio2011}.
Therefore, unlike in the continuous case, $J_h'(\Omega_h; \cdot)$
does not belong, in general, to the image space of $N_h$.
Consequently, 
the solution $\shapedefdisc$
of \eqref{eq:discrete_shape_gradient_1}
has contributions not only from normal forces in the shape derivative $J_h'(\Omega_h;\cdot)$, but also from interior forces as well as tangential boundary forces.
Numerical examples in \cref{sec:numerical_results_gradient} will show that these interior and tangential forces are responsible for spurious descent directions, which in turn lead to degenerate meshes.

Therefore,
we conclude that it is not reasonable to try to solve
\begin{equation*}
	\text{Minimize } J_h(\Omega_h)
\end{equation*}
or its stationarity condition
\begin{equation}
	\label{eq:discrete_stationary_point}
	\text{Find a triangulation }
	\Omega_h
	\text{ such that }
	\shapegraddisc = - E_h^{-1} J_h'(\Omega_h;\cdot) = 0
\end{equation}
as a discretization of
the continuous problem
\eqref{eq:general_shape_optimization_problem}.

Hence, we consider the discretization of \eqref{eq:continuous_shape_gradient_2}
\begin{equation}
	\label{eq:discrete_shape_gradient_2}
	\begin{aligned}
		\text{Minimize}\quad &
		J_h'(\Omega_h; \shapedefdisc)
		+
		\frac12 \, \dual{E_h \shapedefdisc}{\shapedefdisc}
		\\
		\text{with respect to}\quad&
		\shapedefdisc \in S^1(\Omega_h)^d,
		F_h \in S^1(\partial\Omega_h)
		\\
		\text{such that}\quad&
		E_h \, \shapedefdisc - N_h \, F_h = 0
		.
	\end{aligned}
\end{equation}
in which we restrict $E_h \, \shapedefdisc$
to the image space of the discrete normal force operator $N_h$.
As in the continuous setting,
this problem is equivalent to the solution of
\begin{equation}
	\label{eq:continuous_saddle_point_system_discrete}
	\begin{pmatrix}
		E_h & 0 & E_h \\
		0 & 0 & -N_h\adjoint\\
		E_h & -N_h & 0
	\end{pmatrix}
	\begin{pmatrix}
		\shapedefdisc \\ F_h \\ \Pi_h
	\end{pmatrix}
	=
	\begin{pmatrix}
		-J_h'(\Omega_h; \cdot) \\ 0 \\ 0
	\end{pmatrix}
	.
\end{equation}
It is clear that \eqref{eq:continuous_saddle_point_system_discrete} can also be reduced as in \eqref{eq:continuous_saddle_point_system_reduced}.
For later reference,
we mention that the solution
$(\shapeprojgraddisc, F_h, \Pi_h)$ of \eqref{eq:continuous_saddle_point_system_discrete} satisfies
\begin{align}
	\nonumber
	\dual{E_h \, \shapeprojgraddisc}{\shapeprojgraddisc}
	&=
	-\dual{E_h \, \shapeprojgraddisc}{\Pi_h}
	-J_h'(\Omega_h; \shapeprojgraddisc)
	\\&=
	-\dual{N_h \, F_h}{\Pi_h}
	-J_h'(\Omega_h; \shapeprojgraddisc)
	\nonumber
	\\&=
	\label{eq:energy_gradient_directional_derivative}
	-J_h'(\Omega_h; \shapeprojgraddisc)
\end{align}
since $N_h\adjoint \, \Pi_h = 0$ holds.
This shows that $\shapeprojgraddisc$ is always a descent direction for the discrete objective $J_h(\Omega_h;\cdot)$.

As we have seen in \eqref{eq:continuous_shape_gradient_3}
for the continuous setting,
the solution $\shapedefdisc$ of \eqref{eq:discrete_shape_gradient_2}
also solves
\begin{equation}
	\label{eq:discrete_shape_gradient_3}
	\begin{aligned}
		\text{Minimize}\quad &
		\frac12 \, \bigdual{E_h (\shapedefdisc - \shapegraddisc)}{\shapedefdisc - \shapegraddisc}
		\\
		\text{with respect to}\quad&
		\shapedefdisc \in S^1(\Omega_h)^d,
		F_h \in S^1(\partial\Omega_h)
		\\
		\text{such that}\quad&
		E_h \, \shapedefdisc - N_h \, F_h = 0
		,
	\end{aligned}
\end{equation}
where
$\shapegraddisc = -E_h^{-1} J_h'(\Omega_h;\cdot)$
is the solution of 
\eqref{eq:discrete_shape_gradient_1}.
Again,
the solution $\shapeprojgraddisc$ of \eqref{eq:discrete_shape_gradient_3}
can be interpreted as the
projection (w.r.t.\ the $E_h$ inner product)
of $\shapegraddisc$
onto the image space of $E_h^{-1} N_h$.
Therefore,
the notation $\shapeprojgraddisc$ for the solution of 
\eqref{eq:discrete_shape_gradient_2}
is justified.

Our main idea is now to propose, instead of \eqref{eq:discrete_stationary_point},
\begin{equation}
	\label{eq:discrete_stationary_point_2}
	\text{Find a triangulation }
	\Omega_h
	\text{ such that }
	\shapeprojgraddisc = 0
\end{equation}
as an appropriate discrete version of \eqref{eq:general_shape_optimization_problem}.
Note that this is fundamentally different from the ad-hoc discretization \eqref{eq:discrete_stationary_point}
since we neglect the contributions of $J_h'(\Omega_h;\cdot)$
which do not belong to the image space of $N_h$.
We will see via numerical examples
that this problem \eqref{eq:discrete_stationary_point_2}
can be solved
to high accuracy
by an iterative algorithm
using the solution $\shapeprojgraddisc$ of \eqref{eq:discrete_shape_gradient_2}
for the displacement of the triangulation $\Omega_h$
(together with a line search).

For later use, we are going to characterize
stationarity of $\Omega_h$ in the sense of \eqref{eq:discrete_stationary_point_2}.
The deformation $\shapedefdisc = 0$ solves the projection problem \eqref{eq:discrete_shape_gradient_3}
if and only if
\begin{equation*}
	\dual{E_h \, \shapegraddisc}{E_h^{-1} \, N_h \, F_h} = 0
	\qquad\forall F_h \in S^1(\partial\Omega_h).
\end{equation*}
This, in turn, is equivalent to
\begin{equation}
	\label{eq:tangent_shape_gradient}
	\int_{\partial\Omega_h} F_h \, (\shapegraddisc \cdot \nu) \, \ds = 0
	\qquad\forall F_h \in S^1(\partial\Omega_h).
\end{equation}
This means that $\Omega_h$ is stationary in the sense of \eqref{eq:discrete_stationary_point_2} if and only if the usual shape gradient $\shapegraddisc$
is a tangential vector field on $\Omega_h$ in a discrete sense.

We can now state a restricted gradient algorithm for the solution of \eqref{eq:discrete_stationary_point_2}, where we use $\shapeprojgraddisc$ as the deformation field which provides the search direction in the domain transformation.
It is sufficient to utilize a simple 
a backtracking strategy
to comply with the
Armijo condition
\begin{equation}
	\label{eq:armijo_backtracking}
	J_h\bigh(){ (\id + \alpha \, \shapeprojgraddisc)(\Omega_h) }
	\le
	J_h( \Omega_h )
	+
	\sigma \, \alpha \, J_h'(\Omega_h; \shapeprojgraddisc).
\end{equation}
Here, $\sigma \in (0,1)$ is a parameter.

Since we are using the perturbation of identity approach \eqref{eq:perturbation_of_identity_transformation} instead of a more sophisticated family of domain transformations, we also perform a mesh quality control in order to avoid gradient steps which are too large.
To this end, we check that the conditions
\begin{equation}
	\label{eq:mesh_quality}
	\frac12 \le \det(\id + \alpha \, D\shapeprojgraddisc) \le 2
	,
	\qquad
	\norm{ \alpha \, D\shapeprojgraddisc }_F \le 0.3
\end{equation}
are satisfied in every cell throughout the entire domain.
Here, $\norm{\,\cdot\,}_F$ denotes the Frobenius norm of matrices.
The first condition monitors the change of volume of the cell, while the second 
additionally inhibits large changes of the angles.
Note that this amounts to checking three inequalities
per cell.
Due to \eqref{eq:energy_gradient_directional_derivative},
we use
\begin{equation}
	\label{eq:convergence_criterion}
	\dual{E_h\shapeprojgraddisc}{\shapeprojgraddisc}
	=
	-J_h'(\Omega_h; \shapeprojgraddisc)
	\le \varepsilon_{\textup{tol}}^2
\end{equation}
as a convergence criterion
for some small $\varepsilon_{\textup{tol}} > 0$.
These considerations lead to \cref{alg:restricted_gradient_descent}.
\begin{algorithm2e}[htp]
	\SetAlgoLined
	\KwData{Initial domain $\Omega_h$ \\
		Initial step size $\alpha$,
		convergence tolerance $\varepsilon_{\textup{tol}}$,\\
		line search parameters $\beta \in (0,1)$, $\sigma \in (0,1)$
	}
	\KwResult{Improved domain $\Omega_h$ on which \eqref{eq:discrete_stationary_point_2} holds up to $\varepsilon_{\textup{tol}}$}
\For{$i \leftarrow 1$ \KwTo $\infty$}{
	Solve the discrete state equation \eqref{eq:weak_form_example_PDE_discrete} for $u_h$\;
	Solve the discrete adjoint equation \eqref{eq:adjoint_problem_discrete} for $p_h$\;
	Solve \eqref{eq:discrete_shape_gradient_2} for $\shapeprojgraddisc$ with shape derivative $J'(\Omega_h;\cdot)$ from \eqref{eq:shape_derivative_discrete}\;
	\If{${\dual{E_h\shapeprojgraddisc}{\shapeprojgraddisc}} \le \varepsilon_{\textup{tol}}^2$}{
			STOP, the current iterate $\Omega_h$ is almost stationary for \eqref{eq:discrete_stationary_point_2}\;
		}
		Increase step size $\alpha \leftarrow \alpha / \beta$\;
	\While{\eqref{eq:armijo_backtracking} or \eqref{eq:mesh_quality} is violated}{
			Decrease step size $\alpha \leftarrow \beta \, \alpha$\;
		}
		Transform the domain according to $\Omega_h \leftarrow (\id + \alpha \, \shapeprojgraddisc)(\Omega_h)$\;
	}
	\caption{Restricted gradient method for \eqref{eq:discrete_stationary_point_2}.}
	\label{alg:restricted_gradient_descent}
\end{algorithm2e}

\subsection{Comparison with the approach of \cite{SchulzSiebenbornWelker2015:2}}
\label{subsec:comparison_SSW2016}
In \cite{SchulzSiebenbornWelker2015:2},
the authors propose a different way to
convert
$J_h'(\Omega_h; \cdot)$
into
a deformation field from $S^1(\Omega_h)^d$.
Instead of solving \eqref{eq:discrete_shape_gradient_1}
directly, i.e.,
\begin{equation*}
	\text{Find } \shapegraddisc \in S^1(\Omega_h)^d
	\quad\text{s.t.}\quad
	E_h \, \shapegraddisc = -J_h'(\Omega_h, W_h) \quad\forall W_h \in S^1(\Omega_h)^d
	,
\end{equation*}
they
propose
that
``%
Only test functions whose
support includes $\Gamma_{\textup{int}}$ are considered [\ldots]%
'', see \cite[p.~2813]{SchulzSiebenbornWelker2015:2}.
In their problem formulation,
$\Gamma_{\textup{int}}$
corresponds
to $\partial\Omega_h$ in our formulation.
We interpret
this as follows.
We denote by
$D_h : S^1(\Omega_h)^d \to S_0^1(\Omega_h)^d$
the projection operator
defined via
\begin{equation*}
	(D_h \, W_h)(x)
	=
	\begin{cases}
		W_h(x) & \text{if $x$ is an interior node of } \Omega_h \\
		0 & \text{if $x$ is a boundary node of } \Omega_h
	\end{cases}
\end{equation*}
for all nodes $x$ from the mesh $\Omega_h$.
Note that $D_h$ can be represented by a diagonal matrix
(with entries $0$ and $1$)
in the standard basis of $S^1(\Omega_h)^d$.
Then,
the deformation $\shapedefssw$
is computed
via the solution of
\begin{equation*}
	\text{Find } \shapedefssw \in S^1(\Omega_h)^d
	\quad\text{s.t.}\quad
	E_h \, \shapedefssw = -J_h'(\Omega_h, D_h \, W_h) \quad\forall W_h \in S^1(\Omega_h)^d
	.
\end{equation*}
We compare this suggestion with our approach.
\begin{itemize}
	\item
		The deformation $\shapedefssw$ can be computed
		faster than $\shapeprojgraddisc$,
		since the linear system is smaller than
		\eqref{eq:continuous_saddle_point_system_discrete}.

	\item
		The deformation $\shapedefssw$
		cannot be understood as a
		(negative) gradient direction
		generated by an inner product
		on (a subspace of) $S^1(\Omega_h)^d$.
		What is more,
		$\shapedefssw$
		may fail to be a descent direction,
		since
		\begin{equation*}
			J_h'(\Omega_h, \shapedefssw)
			=
			-J_h'\bigl(\Omega_h, E_h^{-1} J_h'(\Omega_h, D_h \, \cdot)\bigr)
		\end{equation*}
		might be positive.
		This indeed does happen in our numerical examples,
		see \cref{sec:numerical_results_gradient} below.
		In contrast,
		our suggestion $\shapeprojgraddisc$
		is generated by an inner product
		on the subspace
		$\{ \shapedefdisc \in S^1(\Omega_h)^d : \exists F_h \in S^1(\partial\Omega_h) : E_h \, \shapedefdisc = N_h \, F_h\}$.

	\item
		The deformation $\shapedefssw$
		is induced by the forces corresponding to
		the linear map
		$W_h \mapsto -J_h'(\Omega_h, D_h \, W_h)$.
		Due to the operator $D_h$,
		these forces only act on the boundary $\partial\Omega_h$.
		However,
		these forces may contain tangential components
		and this is a crucial difference to our approach.
		Eventually,
		this will generate tangential movement of boundary points
		leading to a deterioration
		of the mesh quality.
		We see the onset of this in \cref{fig:gradient_vs_restricted} (center bottom), where boundary nodes accumulate on the right part of the boundary and a rarefaction of nodes occurs on the left.
\end{itemize}

\section{Numerical Results: Comparison of Gradient Methods}
\label{sec:numerical_results_gradient}

The main goal of this section is to compare our proposed restricted gradient method, see \cref{alg:restricted_gradient_descent}, to a classical shape gradient method.
The latter is identical to \cref{alg:restricted_gradient_descent} except that $\shapeprojgraddisc$ is replaced everywhere by the negative shape gradient $\shapegraddisc$ from \eqref{eq:discrete_shape_gradient_1}.
In addition, we also compare it to a method utilizing the deformation fields $\shapedefssw$ obtained along the lines of \cite{SchulzSiebenbornWelker2015:2}.
We refer to the latter as a gradient-like method since it may fail to produce descent directions.
Consequently, $\shapedefssw$ cannot be a negative gradient direction w.r.t.\ any inner product in this situation.

We consider our model problem \eqref{eq:shape_optimization_problem_continuous} with data $f$ as in \cref{example:spurious_descent_direction}.
The line search parameters $\beta= 0.5$ and $\sigma = 0.1$ are used and the initial step size is chosen as $\alpha = 1$.
For the Lam{\'e} and damping parameters in the elasticity operator \eqref{eq:elasticity_operator} we choose
\begin{equation*}
	\mu = \frac{E_0}{2 \, (1+\nu)},
	\quad
	\lambda = \frac{E_0 \, \nu}{(1+\nu)(1-2 \nu)},
	\quad
	\delta = 0.2 \, E_0
\end{equation*}
where $E_0 = 1.0$ is Young's modulus and $\nu = 0.4$ is the Poisson ratio.
The initial shape for all three methods is the same as in \cref{fig:destroy_mesh}~(left).
For this first result, the mesh has 864~triangles and 469~vertices but computations on refined meshes are reported below in \cref{tab:mesh_level}.

We implemented the restricted gradient method, \cref{alg:restricted_gradient_descent}, its classical counterpart as well as the gradient-like method from~\cite{SchulzSiebenbornWelker2015:2} in \fenics, version~2018.1 (\cite{LoggMardalWells2012:1}).
We report computational results obtained on a machine with an Intel(R) Xeon(R) CPU E5--4640 at 2.4~GHz.

Our implementation is freely available on GitHub, see \cite{EtlingHerzogLoayzaWachsmuth2018:2}.
All derivatives were automatically generated by the built-in algorithmic differentiation capabilities of \fenics. 
The restricted shape gradient $\shapeprojgraddisc$, i.e., the solution of \eqref{eq:discrete_shape_gradient_2},
was computed via the discrete counterpart of \eqref{eq:continuous_saddle_point_system_reduced}.
The linear system was solved using \scipy's \lstinline!spsolve! with the \superlu\ solver (\cite{Li2005}), i.e., with the setting \lstinline!use_umfpack = False!.

The restricted gradient method reached the desired tolerance 
\begin{equation}
	\label{eq:restricted_gradient_stopping_criterion}
	\norm{\shapeprojgraddisc}_{E_h}
		=
		\sqrt{\abs{J_h'(\Omega_h; \shapeprojgraddisc)}}
	\le 
	\varepsilon_{\textup{tol}} = 10^{-7}
\end{equation}
at iteration~864 after 40~seconds , while the classical gradient method was stopped at iteration~1500 after 43~seconds , where it had only reached
\begin{equation*}
	\norm{\shapegraddisc}_{E_h}
		=
		\sqrt{\abs{J_h'(\Omega_h; \shapegraddisc)}}
	\approx
	4\cdot 10^{-3}
	.
\end{equation*}
On the other hand, the gradient-like method from \cite{SchulzSiebenbornWelker2015:2} was stopped at iteration 1435 since it failed to generate a descent direction. 
This process took 49~seconds, and at that point the method had reached
\begin{equation*}
		\sqrt{\abs{J_h'(\Omega_h; \shapedefssw)}}
		\approx
		3\cdot 10^{-5}
		.
\end{equation*}

\begin{figure}[ht]
	\begin{subfigure}[b]{\linewidth}
		\centering
		\hspace*{-5mm}
		\includegraphics[width=0.48\linewidth]{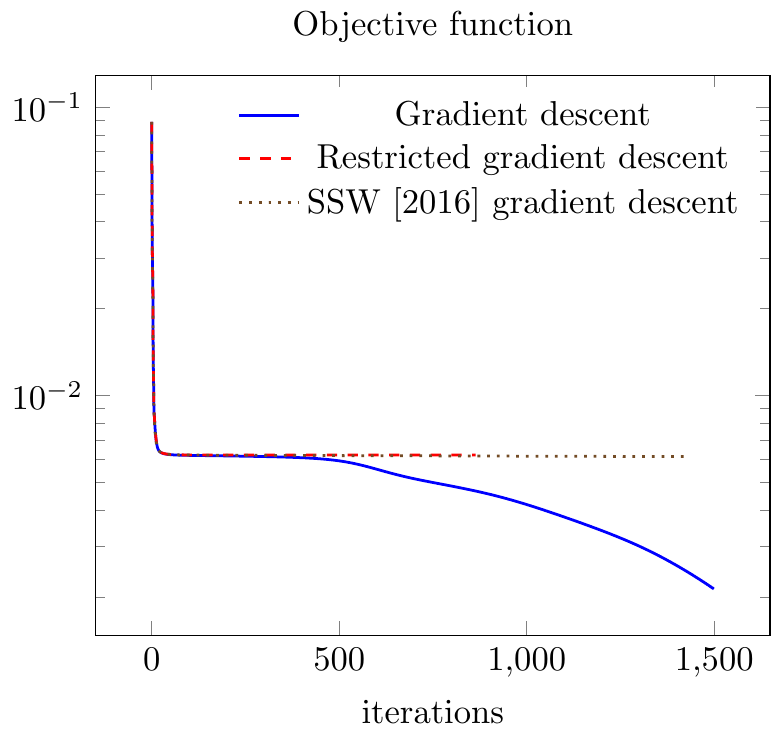}
		\hfill
		\includegraphics[width=0.48\linewidth]{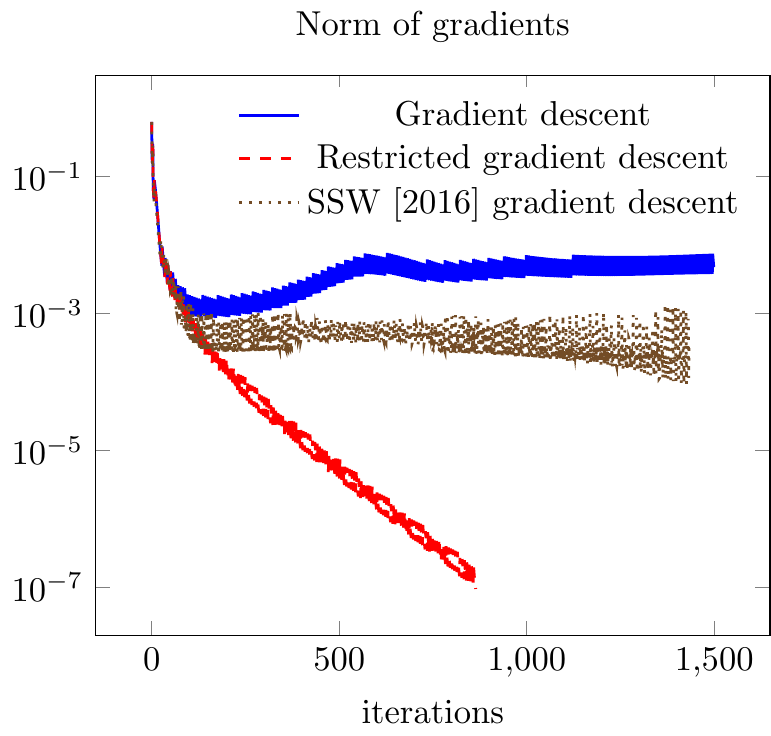}
	\end{subfigure}%
	\caption{History of the objective value $J_h(\Omega_h) + 0.1$ (left) and the norm of the gradients $\norm{\shapegraddisc}_{E_h}$ (for the gradient descent method), $\norm{\shapeprojgraddisc}_{E_h}$ (for the restricted gradient method)
		and $\sqrt{J_h'(\Omega_h; \shapedefssw)}$ (for the gradient-like method from \cite{SchulzSiebenbornWelker2015:2})
	along the iterations (right).}
	\label{fig:history_norm_grad}
\end{figure}

\Cref{fig:history_norm_grad} shows the complete history of the objective and respective shape gradient norms.
The geometry condition \eqref{eq:mesh_quality} was violated only once for all methods, namely in the first very iteration, leading to a reduction of the initial step size. The Armijo condition \eqref{eq:armijo_backtracking} failed approximately once per iteration on average. Typical accepted step sizes were $\alpha = 0.5$ and occasionally $\alpha = 1$.

\begin{figure}[htp]
	\begin{subfigure}[b]{\linewidth}
		\centering
		\includegraphics[width=.33\textwidth]{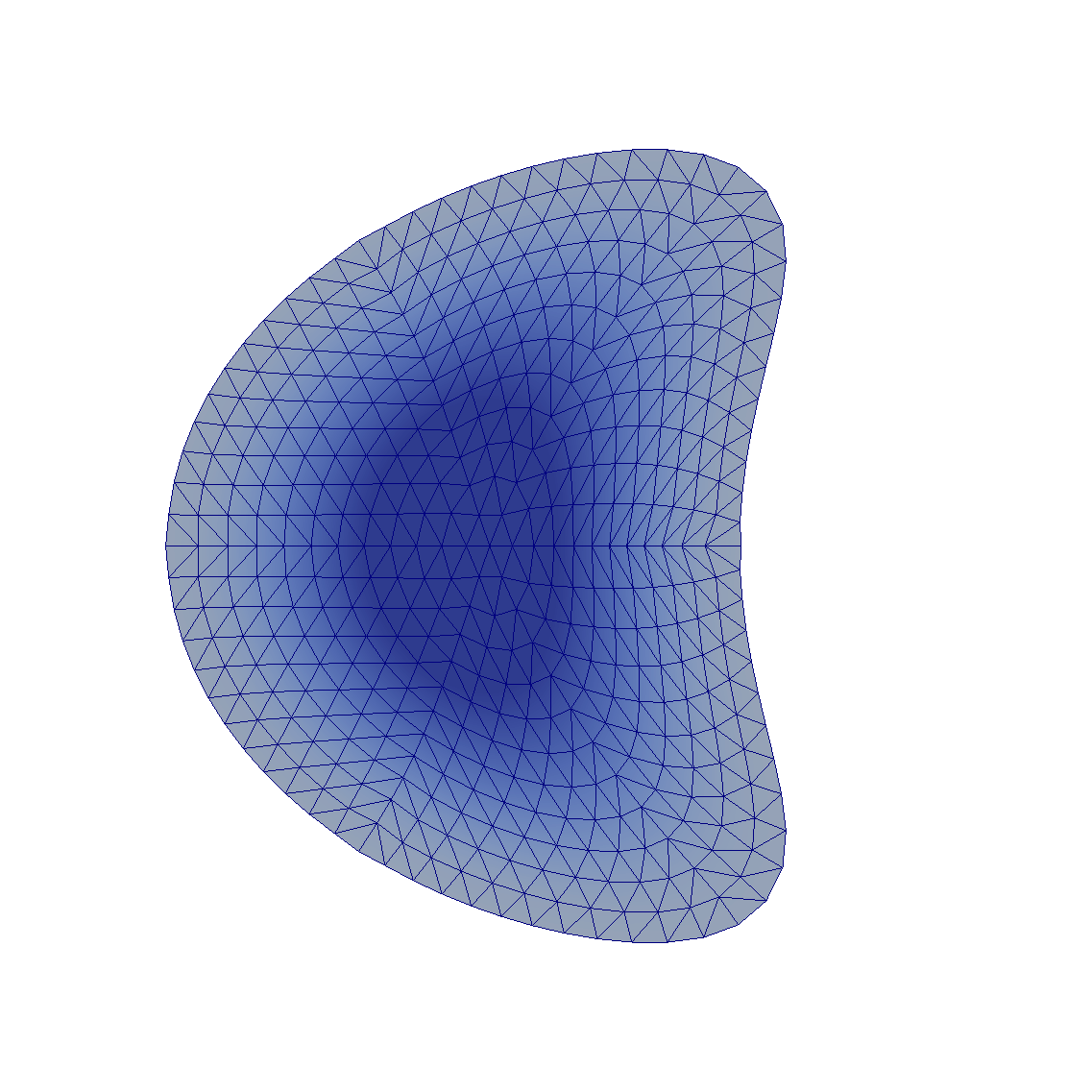}%
		\hfill%
		\includegraphics[width=.33\textwidth]{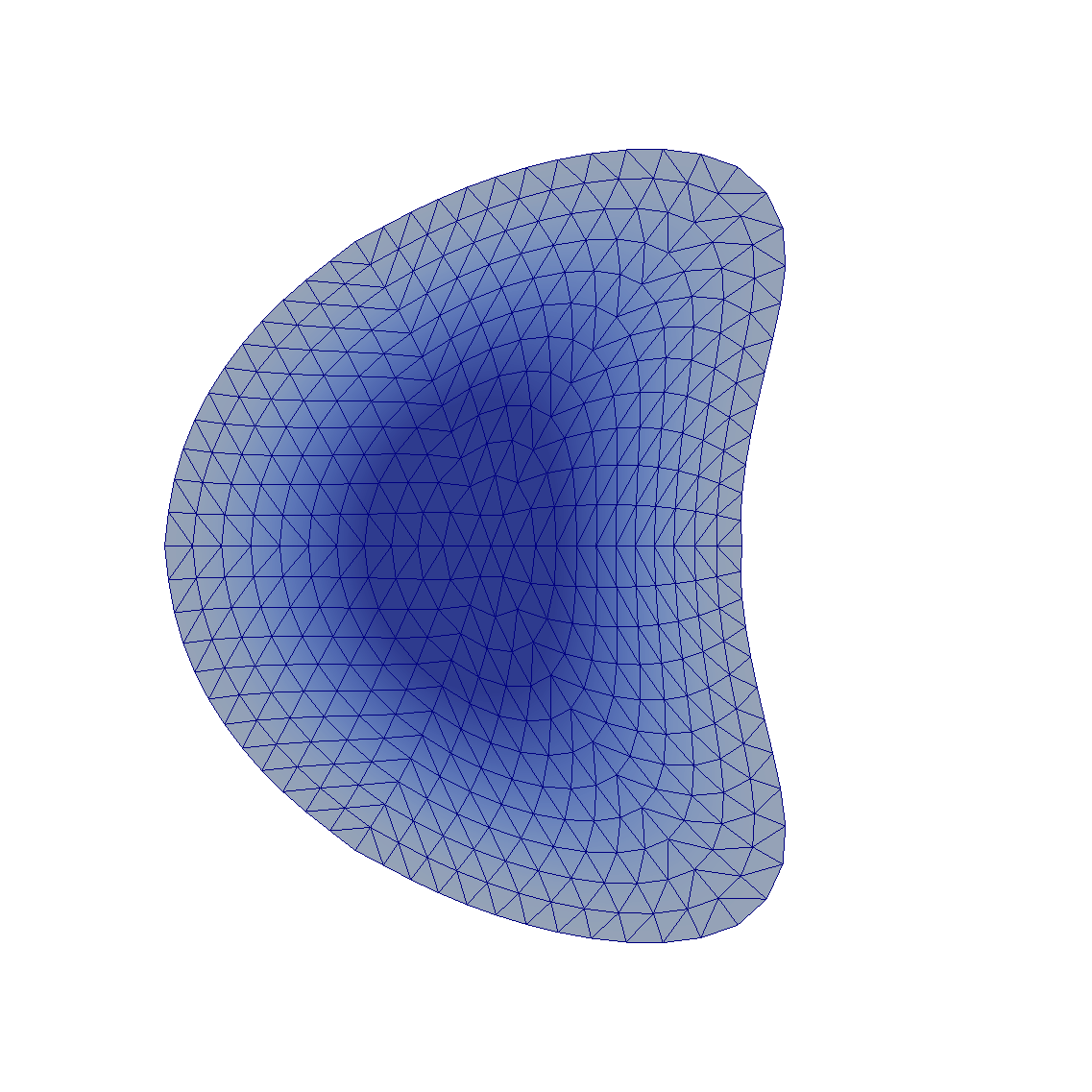}%
		\hfill%
		\includegraphics[width=.33\textwidth]{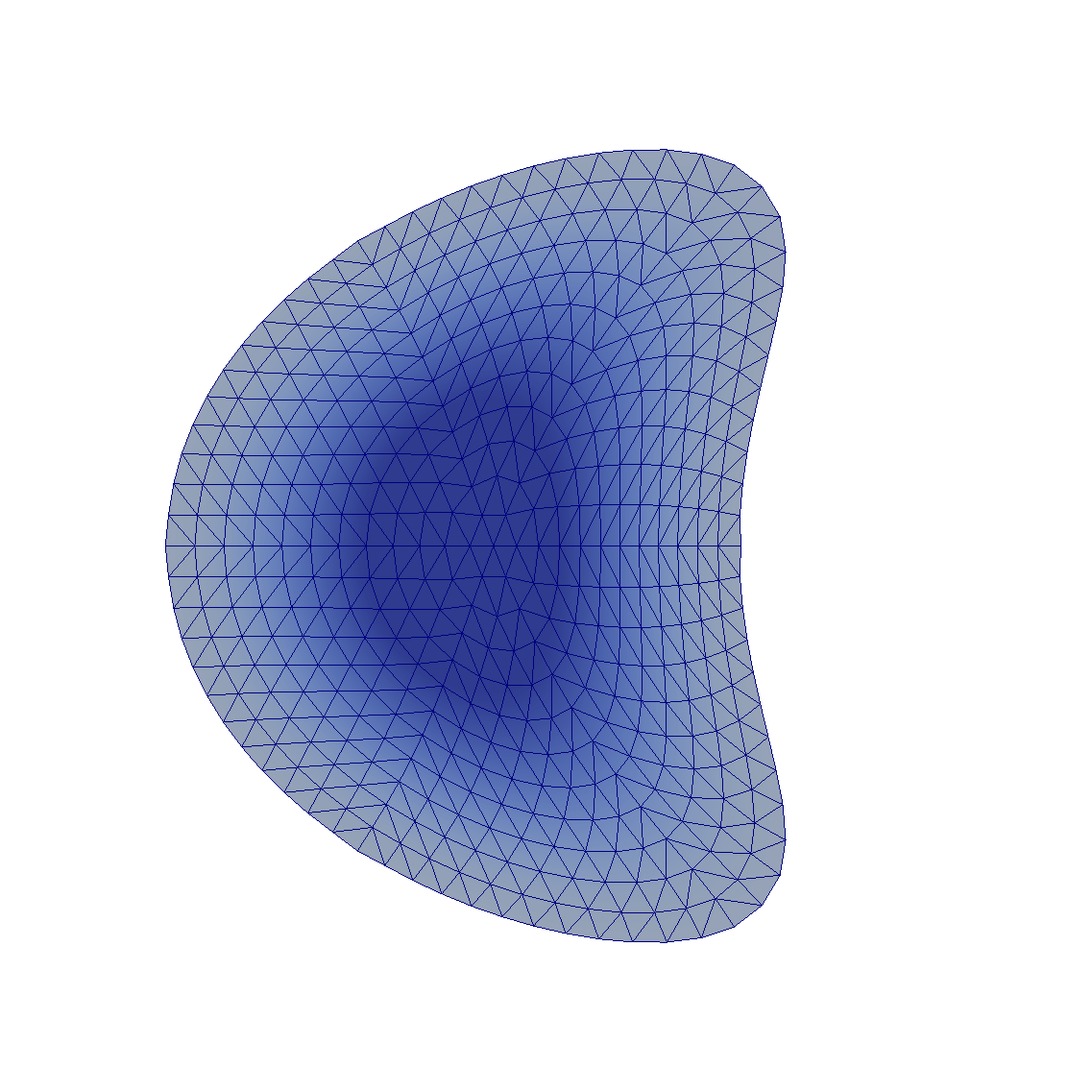}%
	\end{subfigure}%
	\\[-.5cm]
	\begin{subfigure}[b]{\linewidth}
		\centering
		\includegraphics[width=.33\textwidth]{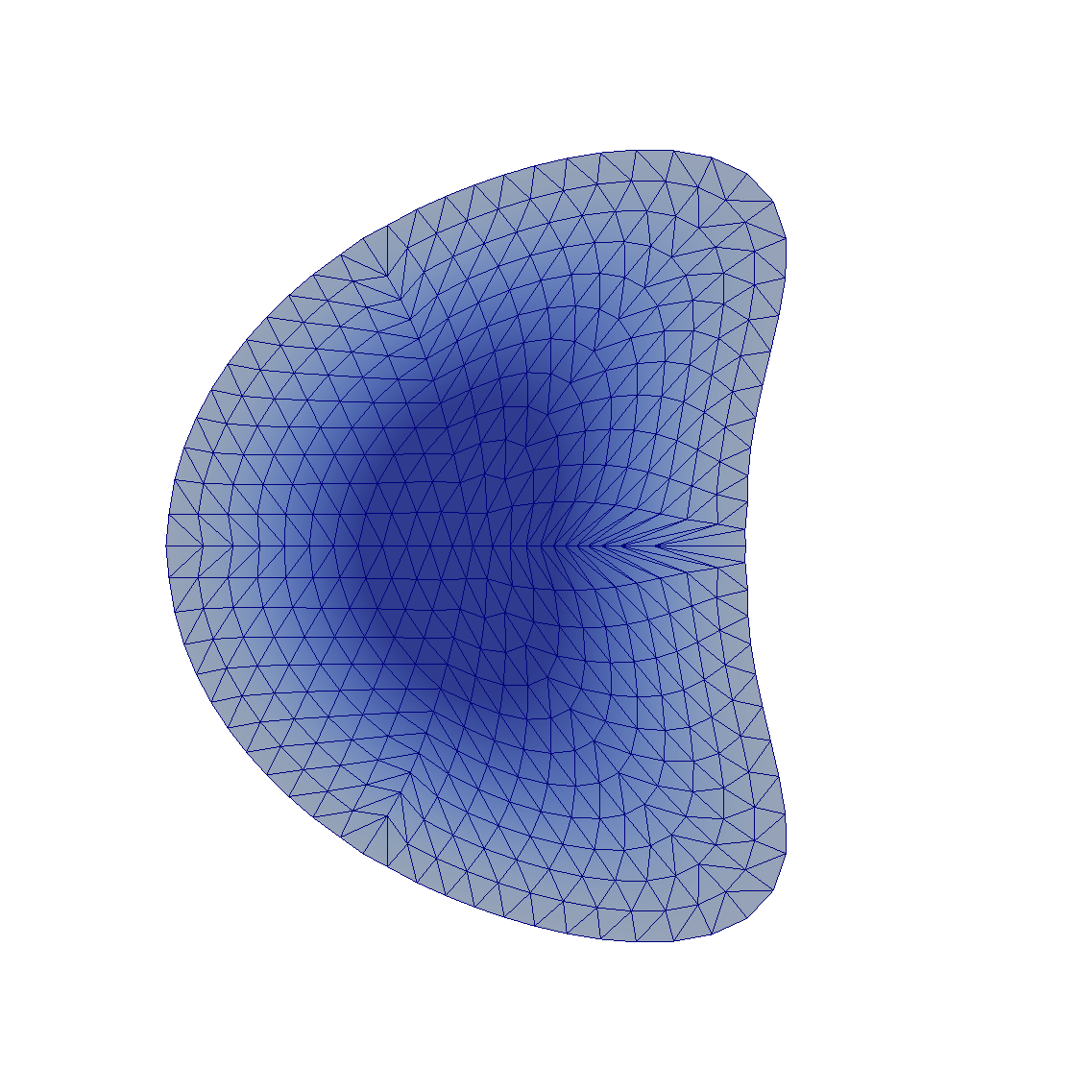}%
		\hfill%
		\includegraphics[width=.33\textwidth]{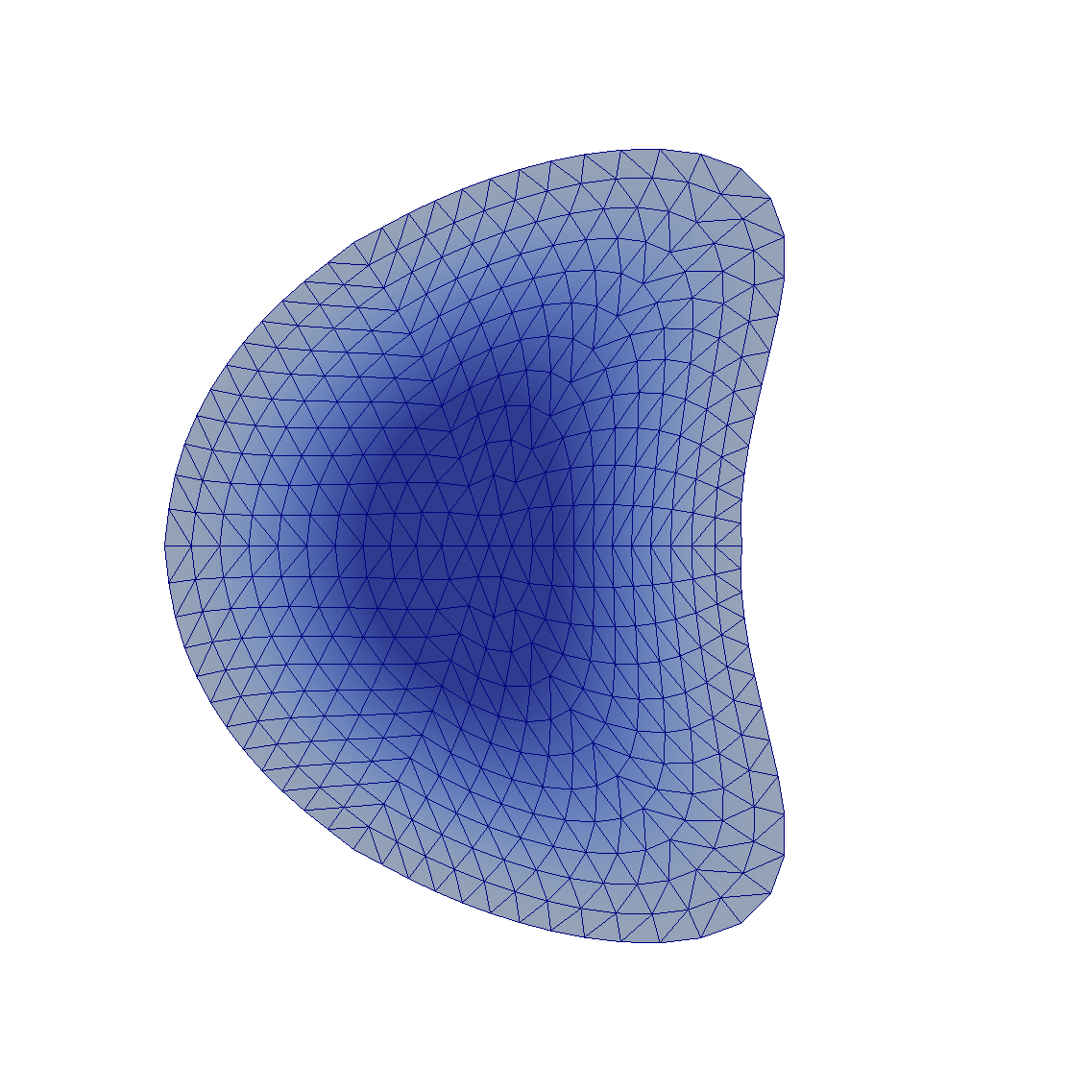}%
		\hfill%
		\includegraphics[width=.33\textwidth]{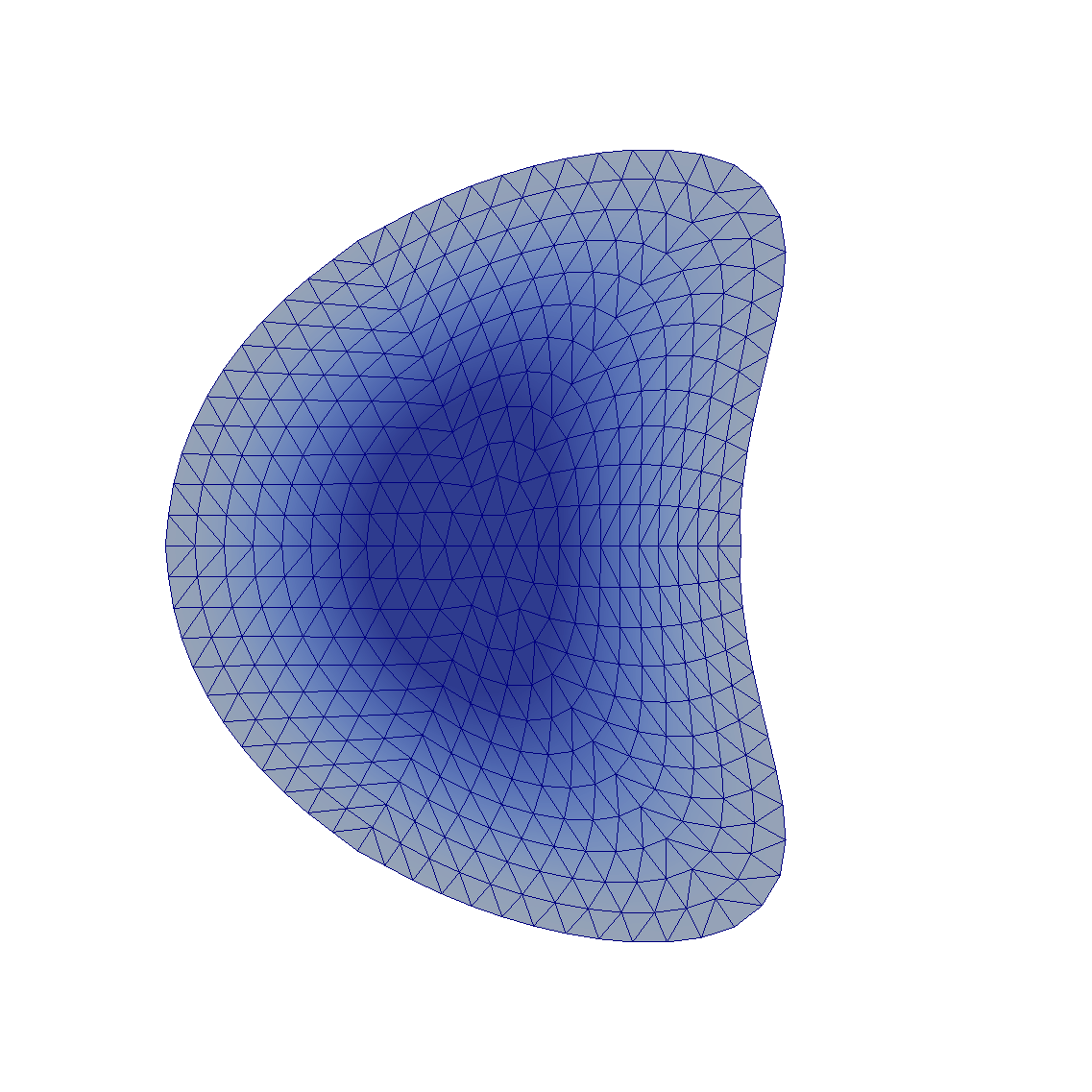}%
	\end{subfigure}%
	\\[-.5cm]
	\begin{subfigure}[b]{\linewidth}
		\centering
		\includegraphics[width=.33\textwidth]{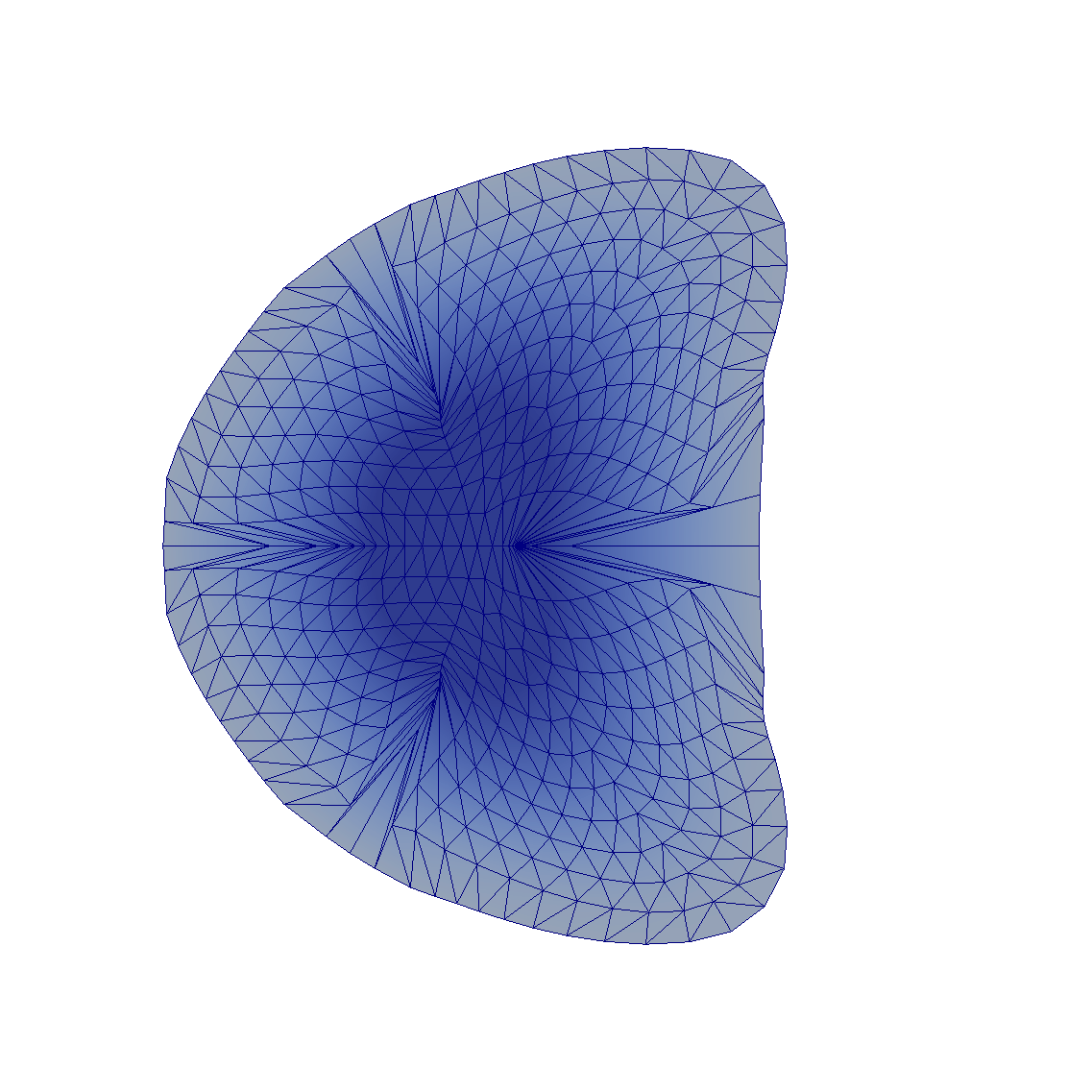}%
		\hfill%
		\includegraphics[width=.33\textwidth]{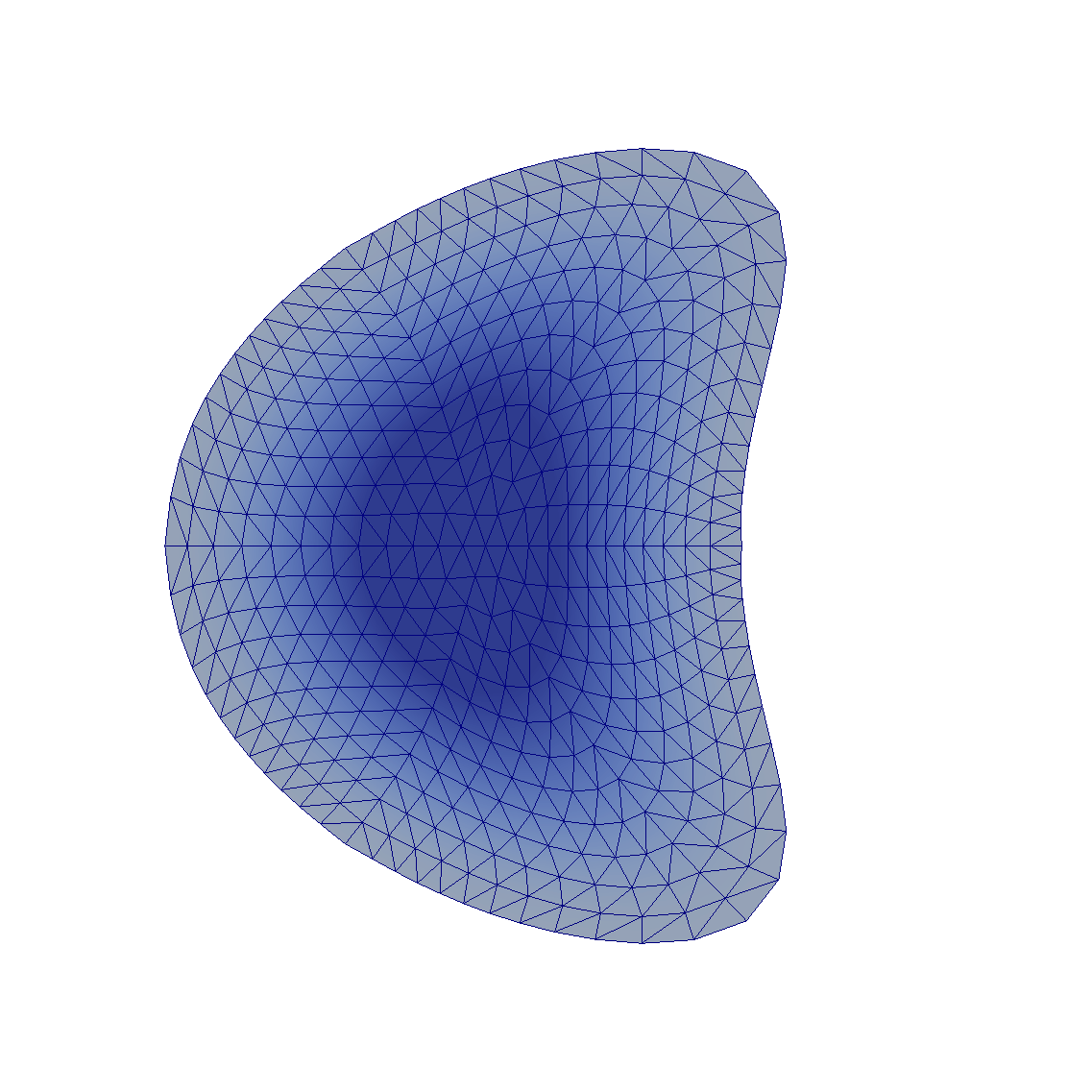}%
		\hfill%
		\includegraphics[width=.33\textwidth]{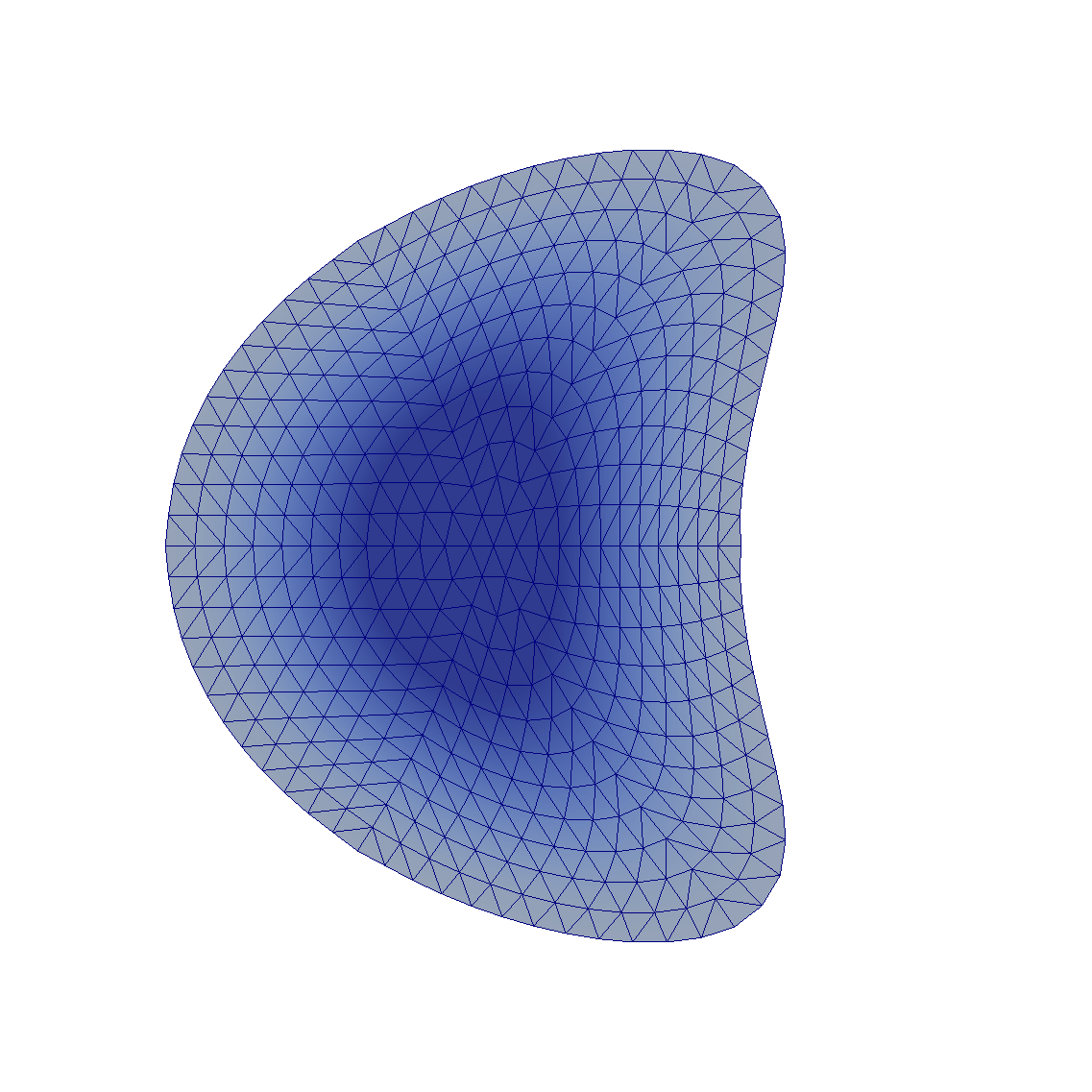}%
	\end{subfigure}
	\caption{Intermediate shapes $\Omega_h$ obtained with the classical gradient method (left), the gradient method from \cite{SchulzSiebenbornWelker2015:2} (middle) and the restricted gradient method (right) at iterations~300, 600, and 1500, 1435 and 864, respectively.
}
	\label{fig:gradient_vs_restricted}
\end{figure}
\cref{fig:gradient_vs_restricted} shows the domains $\Omega_h$ during the iteration of all three methods for comparison.
It can clearly be inferred that the initial iterates are virtually identical.
The classical gradient method begins to produce visibly different shapes around iteration~500, when the objective value (shown in \cref{fig:history_norm_grad}) has practically converged but the gradient norms are still 
\begin{equation*}
	\norm{\shapegraddisc}_{E_h}
	\approx
	5\cdot 10^{-3}
	\quad
	\text{and}
	\quad 
	\norm{\shapeprojgraddisc}_{E_h}
	\approx
	4 \cdot 10^{-6},
\end{equation*}
respectively.
At this point, the classical gradient method starts to pursue spurious descent directions, which results in a further decrease of the discrete objective at the expense of increasingly degenerate meshes.
Similarly,
the gradient method from \cite{SchulzSiebenbornWelker2015:2} produces some
tangential movement on the boundary.
This decreases the mesh quality slightly
and inhibits further decrease of the norm of the gradient.
\begin{figure}[ht]
	\begin{subfigure}[b]{\linewidth}
		\centering
		\hspace*{-5mm}
		\includegraphics[width=0.48\linewidth]{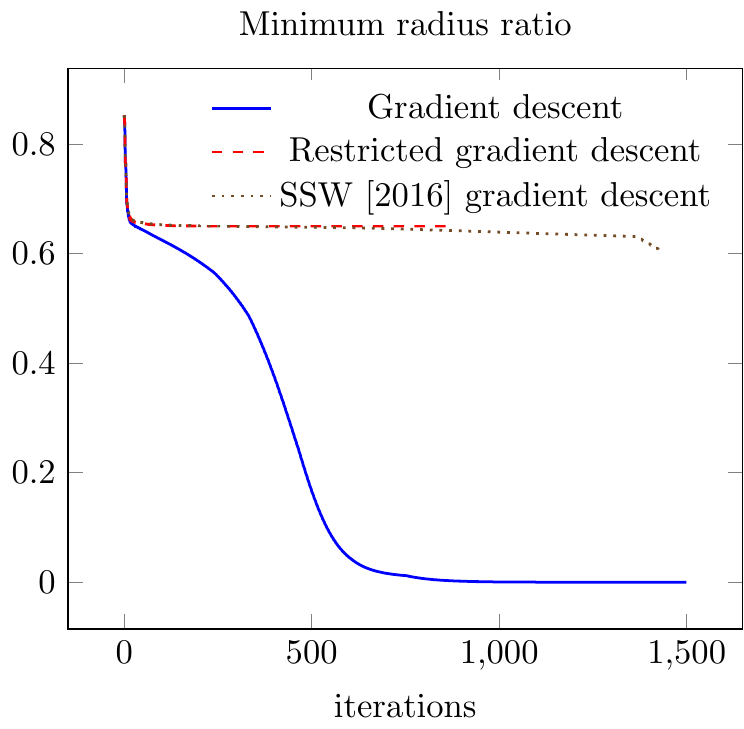}
	\end{subfigure}%
	\caption{History of the mesh quality indicator (minimum over all cells of two times the inradius divided by the circumradius of the cell).}
	\label{fig:history_mesh_quality}
\end{figure}
As an indicator for the mesh quality,
	we used the ``minimum radius ratio'',
	which is provided by \fenics\ and is defined as the minimum over all cells of two (the geometric dimension) times the inradius divided by the circumradius of the cell.
	This value lies between $0$ and $1$, where $0$ corresponds to a degenerate cell
	and $1$ to an equilateral triangle.
	The results for the mesh quality are shown in \cref{fig:history_mesh_quality}.%

	We conclude from this numerical experiment that the restricted gradient method is slightly more expensive per iteration compared to the classical gradient method and the one obtaining its search direction from \cite{SchulzSiebenbornWelker2015:2}.
	However, the proposed restricted gradient method reduces the norm of the restricted gradient much more effectively than the other two methods.
	The cell aspect ratios for the restricted gradient method and the one from \cite{SchulzSiebenbornWelker2015:2} are equally good. 
However, the latter eventually failed to produce descent directions in our experiment, while the classical gradient method created a distorted mesh and did not converge.

\begin{figure}
	\includegraphics[scale=0.11]{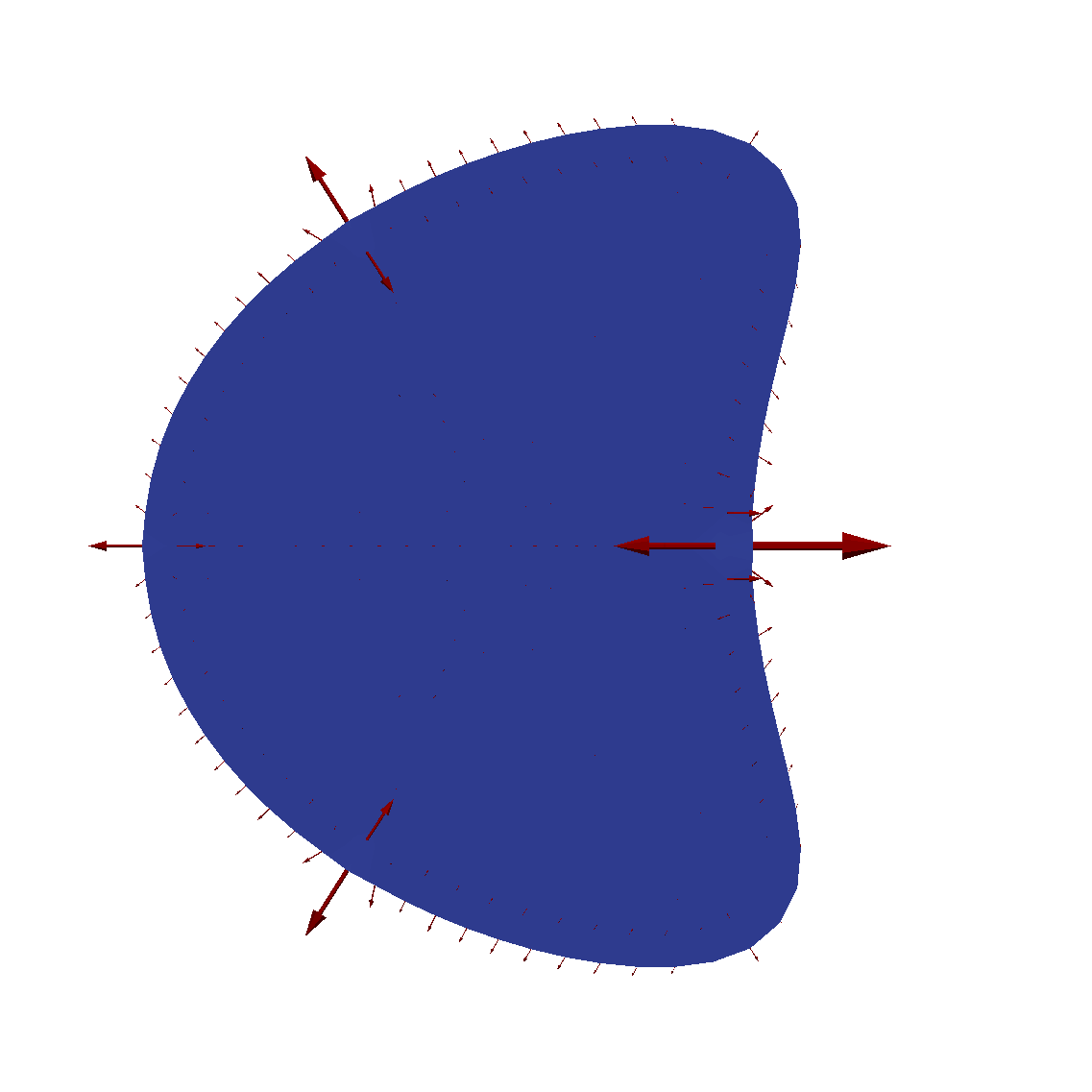}%
	\hfill%
	\includegraphics[scale=0.11]{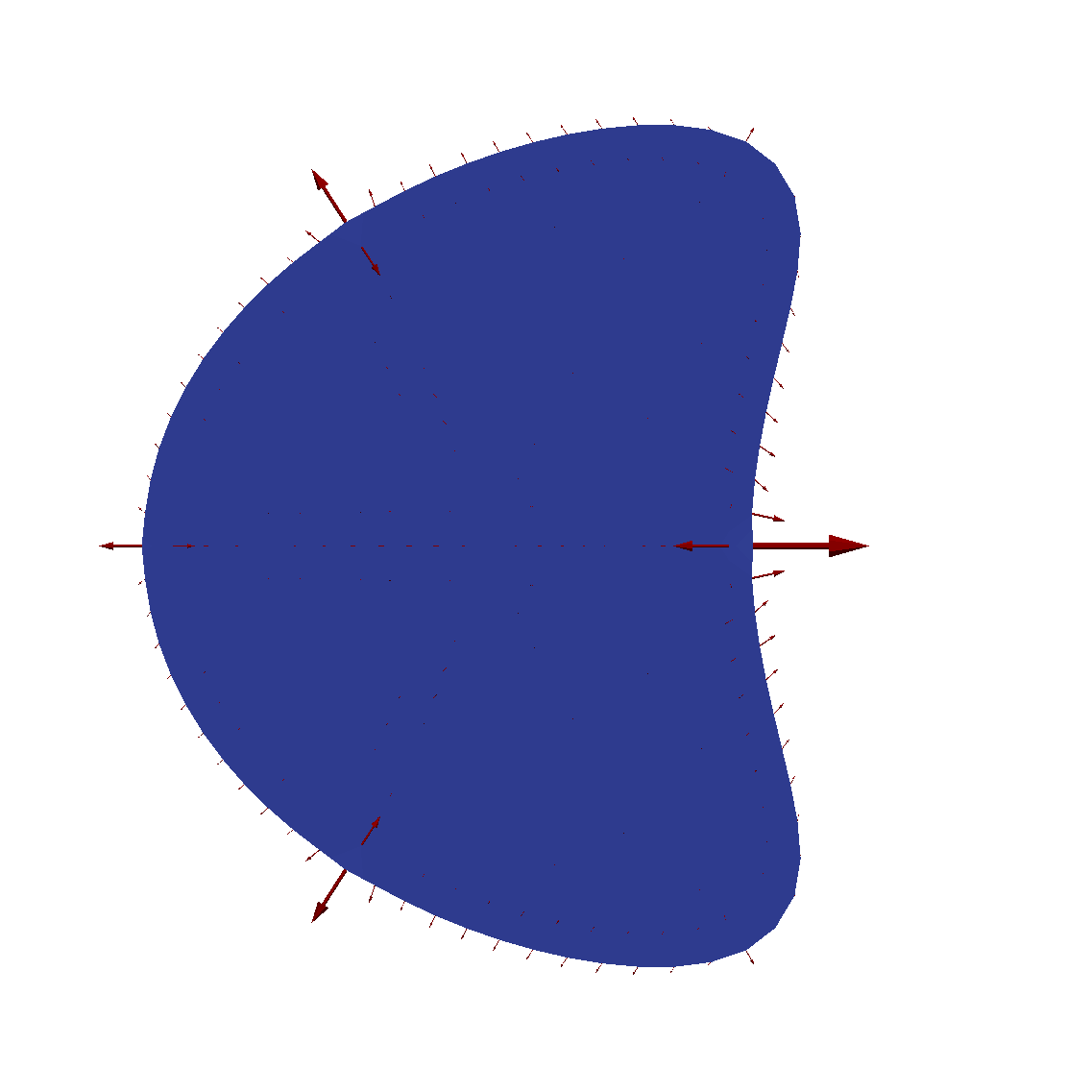}%
	\\%
	\includegraphics[scale=0.11]{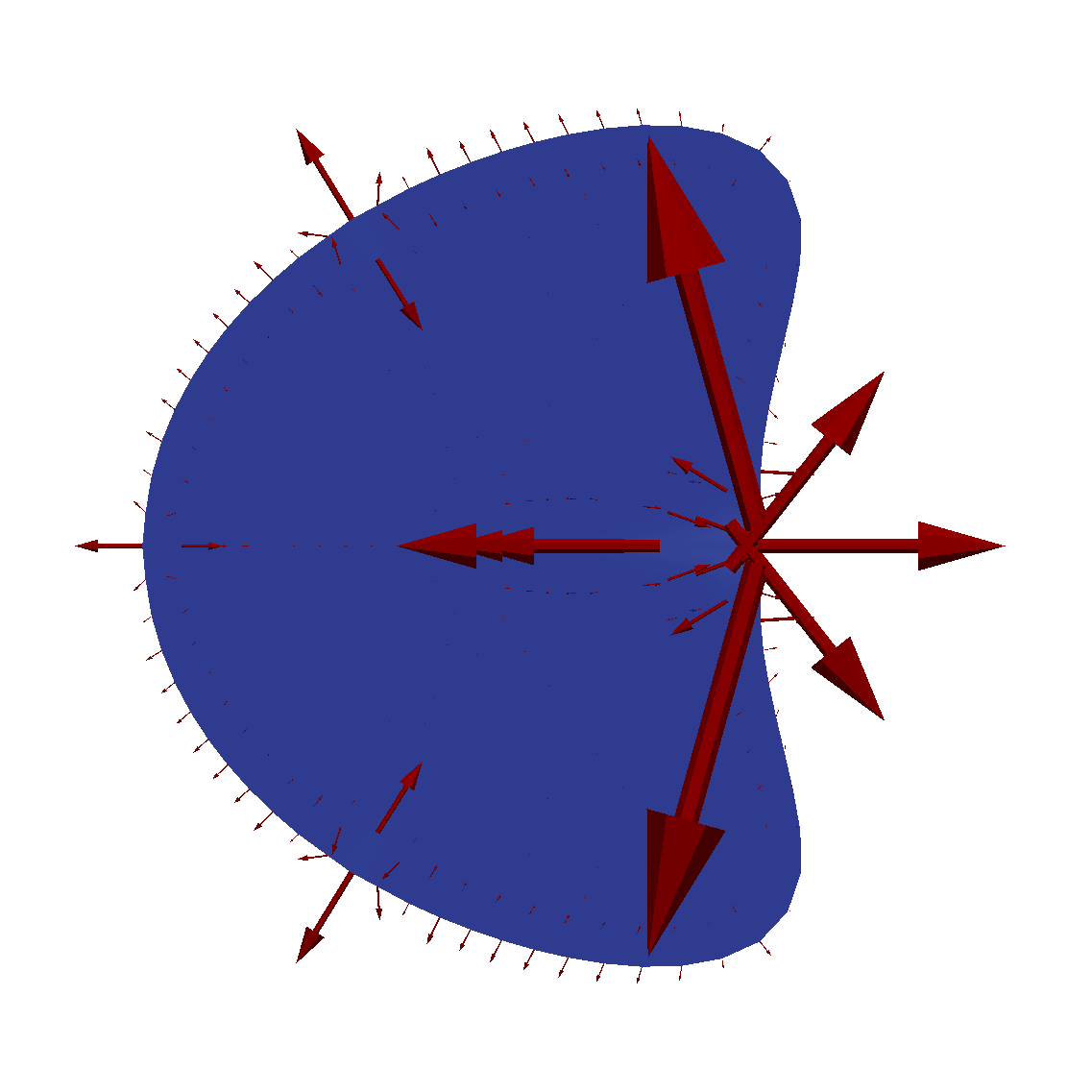}%
	\hfill%
	\includegraphics[scale=0.11]{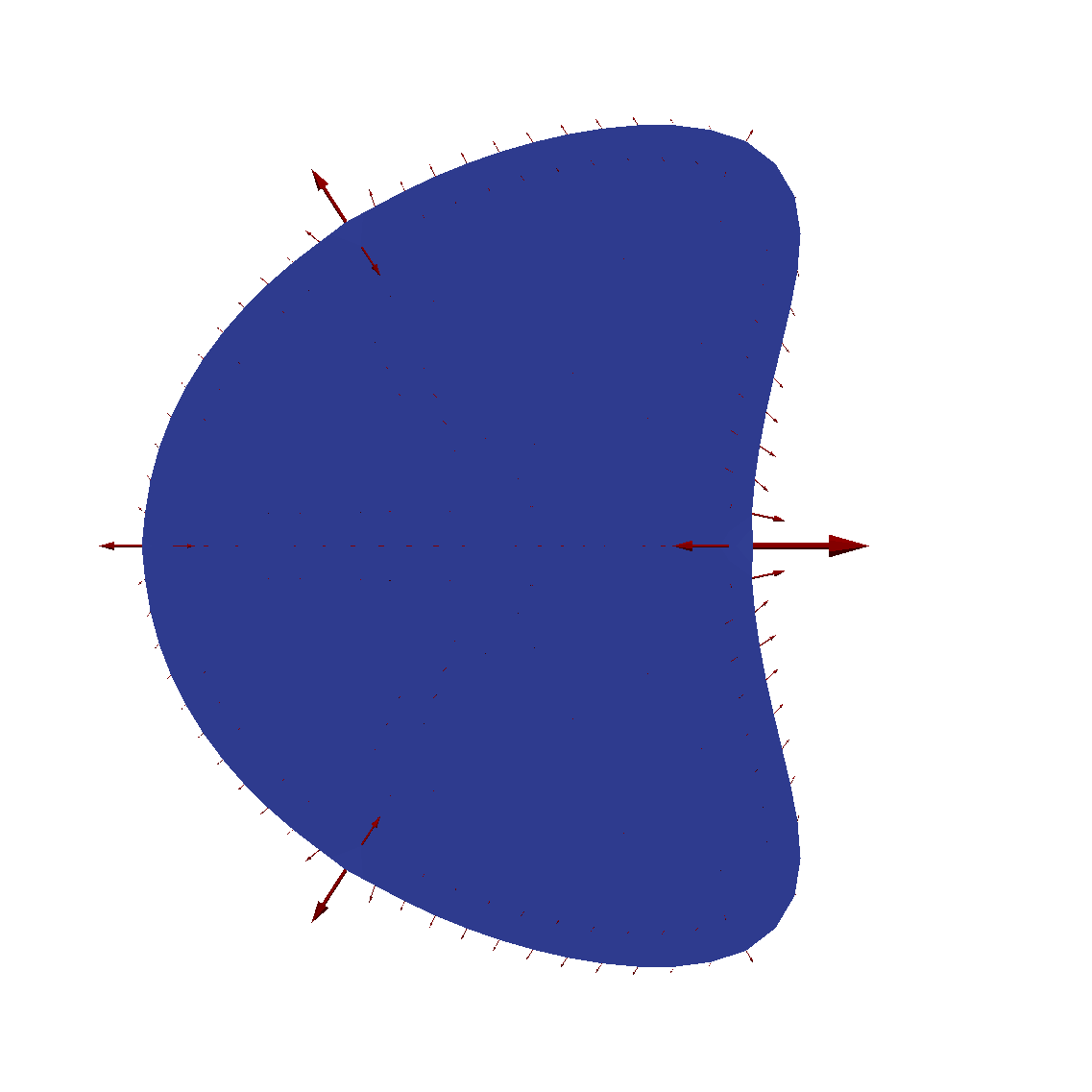}%
	\\%
	\includegraphics[scale=0.11]{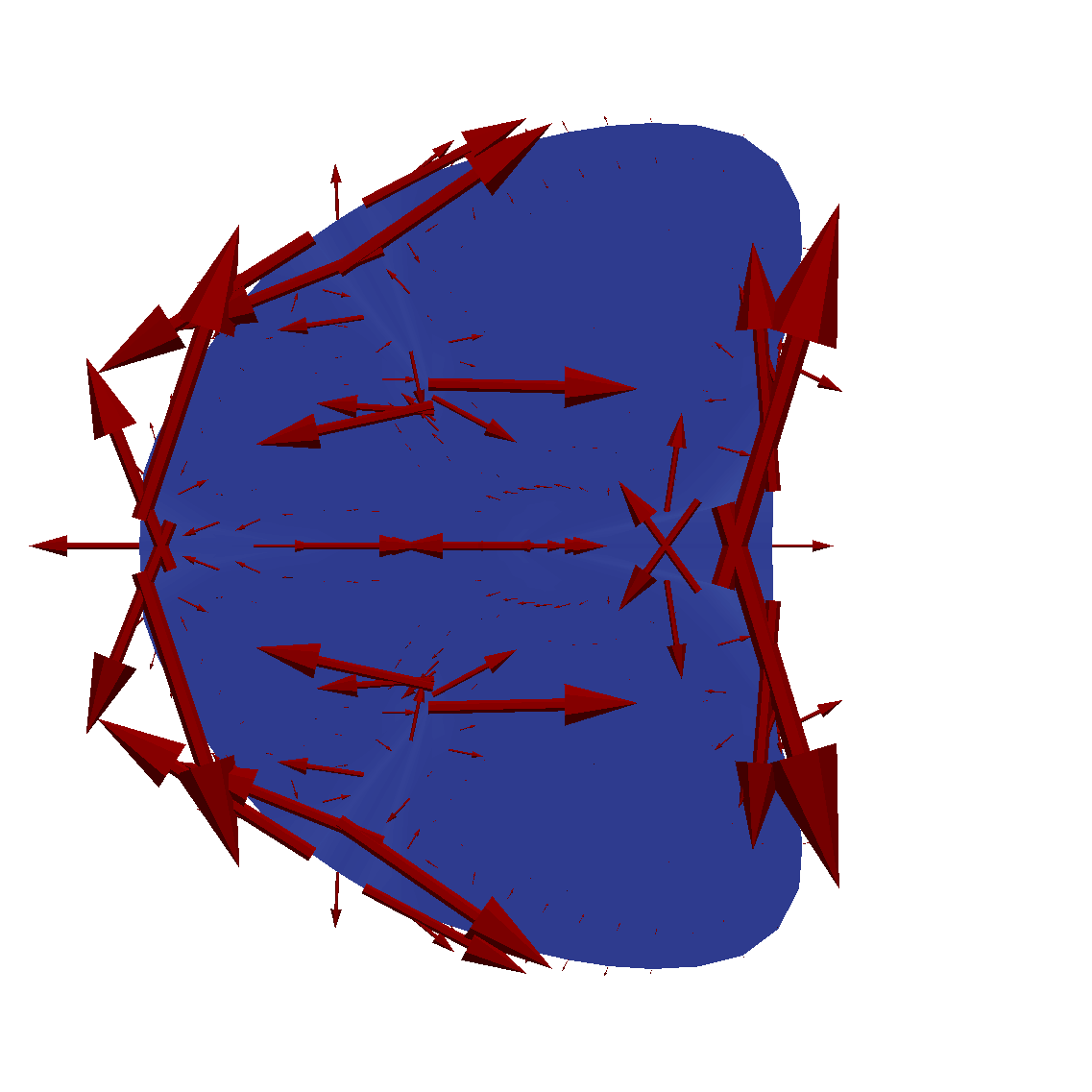}%
	\hfill%
	\includegraphics[scale=0.11]{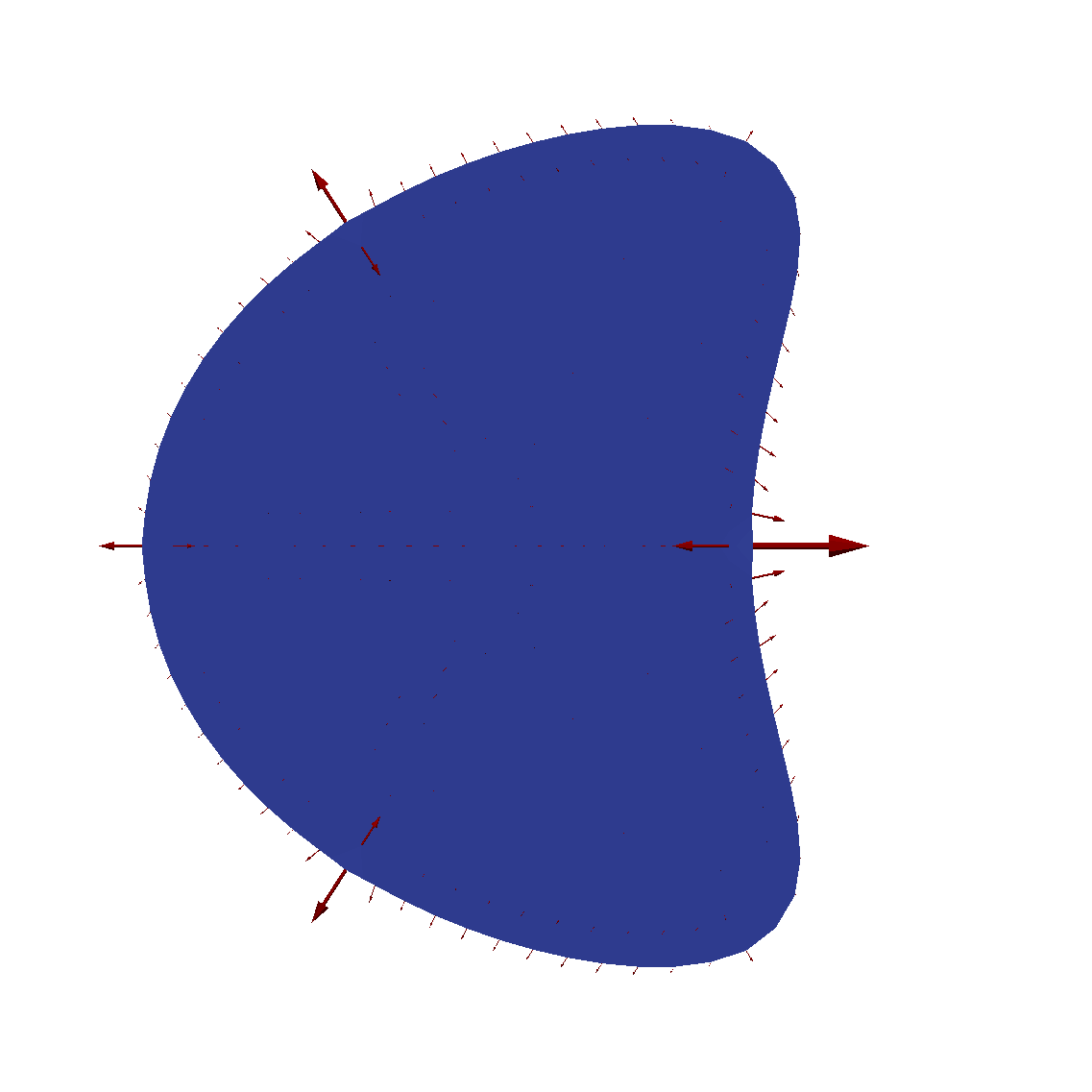}%
	\caption{Visualization of the shape derivatives $J_h'(\Omega_h;\cdot)$ obtained with the classical gradient method (left), and the restricted gradient method (right) at iterations 300, 600, and 1500, and 864, respectively.}
	\label{fig:intermediate_shapesderivatives}
\end{figure}

To further illustrate this point, we show in \cref{fig:intermediate_shapesderivatives} visualizations of the shape derivative $J_h'(\Omega_h;\cdot)$ for the classical and restricted gradient methods; see \eqref{eq:shape_derivative_discrete}.
In fact, this is a linear functional on the space of piecewise linear perturbation fields $\shapedefdisc \in S^1(\Omega_h)^d$.
In \cref{fig:intermediate_shapesderivatives} we display the $S^1(\Omega_h)^d$ representer of $J_h'(\Omega_h;\cdot)$ w.r.t.\ the $L^2$ inner product, i.e., we solve a linear system governed by a block-diagonal mass matrix.

Let us comment on the shape derivative for the restricted gradient method as shown in the right column of \cref{fig:intermediate_shapesderivatives}.
It is apparent that the displacement field $\shapegraddisc$, i.e., the solution of \eqref{eq:discrete_shape_gradient_1}, is non-zero and in fact essentially the same for the iterations~300, 600, and 864 shown.
However $\shapegraddisc$ also has essentially no component in the space of deformations induced by normal forces.
Therefore its projection into this space, see \eqref{eq:discrete_shape_gradient_3}, leaves us with a very small norm $\norm{\shapeprojgraddisc}_{E_h}$, as shown in \cref{fig:history_norm_grad}.
The images visualizing the shape derivative for the classical gradient method in the left column of \cref{fig:intermediate_shapesderivatives} show that the method has allowed the spurious part of the derivative to build up, which eventually dominates the search direction.

\section{Restricted Newton-Like Method}
\label{sec:restricted_shape_hessian}

In the previous two sections
we have seen that \eqref{eq:discrete_stationary_point_2}
is a reasonable discrete optimality condition 
and that it can be solved to high accuracy via a first-order gradient descent method.
However, as is well known for the minimization of even mildly ill-conditioned quadratic polynomials, gradient descent methods require a large number of iterations to achieve convergence.
We observed the same behavior in \cref{sec:numerical_results_gradient}.

Therefore, we are also investigating a Newton-like method
for solving \eqref{eq:discrete_stationary_point_2}.
First, we focus on the continuous case
and comment on its discretization afterwards.
Let $\Omega$ be our current iterate.
As before, we denote by $u$ the associated state, see \eqref{eq:weak_form_example_PDE}, and by $p$ the adjoint state, see \eqref{eq:adjoint_problem}.
The solution of the restricted shape gradient problem \eqref{eq:continuous_saddle_point_system} at $\Omega$ is denoted by
$(\shapeprojgrad, F, \Pi)$.
Recall that our goal is to achieve $\shapeprojgrad = 0$ or, equivalently,
$F = 0$, cf.\ \eqref{eq:discrete_stationary_point_2}.
In practice, we impose a stopping criterion of the form 
$\norm{\shapeprojgrad}_E \le \varepsilon_{\textup{tol}}$
as we did for the gradient method.

In order to allow the reader to follow the derivation for the solution of \eqref{eq:discrete_stationary_point_2} of our Newton method more easily, we draw the parallel with Newton's method for $\Phi(x) = 0$ for some $\Phi: \R^n \to \R^n$.
We consider the equation $\Phi(x + \delta x) = 0$ for the unknown update $\delta x$.
In our context the iterate $x$ represents the current domain $\Omega$ and the update corresponds to a perturbation field $\newtonupdate$.
Since the update takes $\Omega$ into a new domain, we need to manipulate the expression $\Phi(x + \delta x) = 0$ and pull it back to $\Omega$.
Finally, we linearize about $\delta x = 0$, which amounts to $\Phi(x) + D\Phi(x) \, \delta x = 0$.

In our Newton method we seek a deformation field $\newtonupdate$
(taking the role of $\delta x$ above)
such that the updated domain
$\Omega_\newtonupdate := (\id + \newtonupdate)(\Omega)$
is stationary in the sense that the solution of \eqref{eq:continuous_saddle_point_system} (at $\Omega_W$ instead of $\Omega$) satisfies $\shapeprojgradnewton = 0$. 
As in \cref{sec:restricted_mesh_deformations} we are only considering updates $\newtonupdate$
which are induced by a normal force $G$,
i.e., $E\,\newtonupdate - N\,G = 0$ should hold.

In order to characterize the stationarity of
the transformed domain
$\Omega_{\newtonupdate}$,
we introduce the elasticity operator
$E_\newtonupdate : H^1(\Omega_\newtonupdate)^d \to (H^1(\Omega_\newtonupdate)^d)\dualspace$
and the normal force operator
$N_\newtonupdate : L^2(\partial \Omega_\newtonupdate) \to (H^1(\Omega_\newtonupdate)^d)\dualspace$
on $\Omega_\newtonupdate$
analogously to
\eqref{eq:elasticity_operator} and \eqref{eq:normal_force_operator}.
With the transformation field $T_\newtonupdate := \id + \newtonupdate : \Omega \to \Omega_\newtonupdate$,
we define the pull-backs $E^\newtonupdate : H^1(\Omega)^d \to (H^1(\Omega)^d)\dualspace$
of $E_\newtonupdate$
and $N^\newtonupdate : L^2(\partial\Omega) \to (H^1(\Omega)^d)\dualspace$ of $N_\newtonupdate$
via
\begin{align*}
	\dual{E^\newtonupdate \, W_1}{W_2}
	&:=
	\dual{E_\newtonupdate \, (W_1 \circ T_\newtonupdate^{-1})}{W_2 \circ T_\newtonupdate^{-1}}
	,
	\\
	\dual{N^\newtonupdate \, F}{W_2}
	&:=
	\dual{N_\newtonupdate \, (F \circ T_\newtonupdate^{-1})}{W_2 \circ T_\newtonupdate^{-1}}
\end{align*}
for $W_1, W_2, W \in H^1(\Omega)^d$.

Since we wish to achieve conditions defined on the current domain $\Omega$, rather than on the unknown transformed domain $\Omega_\newtonupdate$ after the Newton step, we consider 
the Lagrangian
associated with problem \eqref{eq:shape_optimization_problem_continuous} on $\Omega_\newtonupdate = T_\newtonupdate(\Omega)$ 
and pull it back to $\Omega$.
Using the usual integral substitution and chain rule, and denoting the pulled-back solutions of the state and adjoint equations on $\Omega_\newtonupdate$ by $u^\newtonupdate$ and $p^\newtonupdate$, respectively, we obtain
$\LL : W^{1,\infty}(\Omega)^d \times H_0^1(\Omega) \times H_0^1(\Omega) \to \R$,
defined as
\begin{align*}
	\MoveEqLeft
	\LL(\newtonupdate, u^\newtonupdate, p^\newtonupdate)
	=
	\int_\Omega u^\newtonupdate \, \det(\id + D\newtonupdate) \, \dx
	\\
	&
	+
	\int_\Omega ((\id + D\newtonupdate)^{-\top} \, \nabla u^\newtonupdate) \cdot ((\id + D\newtonupdate)^{-\top} \, \nabla p^\newtonupdate) \, \det(\id + D\newtonupdate) \, \dx
	\\
	&
	-
	\int_\Omega (f \circ T_{\newtonupdate}) \, p^\newtonupdate \, \det(\id + D\newtonupdate) \, \dx
	.
\end{align*}
Notice that $\frac{\partial}{\partial \newtonupdate} \LL( \newtonupdate, u^{\newtonupdate}, p^{\newtonupdate} )$ is the shape derivative $J'(\Omega_\newtonupdate;\cdot)$.
Thus we find that the stationarity of $\Omega_\newtonupdate$
is equivalent to the requirement
that the solution of the nonlinear system
\begin{equation*}
	\begin{pmatrix}
		\cdot & \cdot & \cdot       & \cdot        & \cdot \\
		\cdot & \cdot & \cdot       & \cdot        & \cdot \\
		\cdot & \cdot & E^{\newtonupdate} & \cdot        & E^{\newtonupdate} \\
		\cdot & \cdot & \cdot       & \cdot        & -(N^{\newtonupdate})\adjoint \\
		\cdot & \cdot & E^{\newtonupdate} & -N^{\newtonupdate} & \cdot
	\end{pmatrix}
	\begin{pmatrix}
		u^{\newtonupdate} \\
		p^{\newtonupdate} \\
		\shapeprojgradnewton \\
		F^{\newtonupdate} \\
		\Pi^{\newtonupdate}
	\end{pmatrix}
	+
	\begin{pmatrix}
		\frac{\partial}{\partial u}             \LL( \newtonupdate, u^{\newtonupdate}, p^{\newtonupdate} ) \\
		\frac{\partial}{\partial p}             \LL( \newtonupdate, u^{\newtonupdate}, p^{\newtonupdate} ) \\
		\frac{\partial}{\partial \newtonupdate} \LL( \newtonupdate, u^{\newtonupdate}, p^{\newtonupdate} ) \\
		\cdot \\
		\cdot
	\end{pmatrix}
	=
	0
\end{equation*}
satisfies $\shapeprojgradnewton = 0$.
In view of the injectivity of $N^W$, this is equivalent to $F^{\newtonupdate} = 0$.
Here, ``$\cdot$'' stands for a zero block.
We mention that the first two equations in this system correspond
to the adjoint and state equation on $\Omega_\newtonupdate$ but pulled back to $\Omega$, respectively.
Moreover, note that the solution
$(u^{\newtonupdate}, p^{\newtonupdate}, \shapeprojgradnewton, F^{\newtonupdate}, \Pi^{\newtonupdate})$ of the above system
is the pull-back of the solution
of the state equation, adjoint equation and the projected shape gradient of the system \eqref{eq:continuous_saddle_point_system}
formulated on the domain $\Omega_\newtonupdate$.

Together with the requirement that the deformation field $\newtonupdate$
itself is induced by some normal force $G$, we have to solve the nonlinear system
\begin{equation}
	\label{eq:massive_nonlinear_system_2}
	\begin{pmatrix}
		\cdot & \cdot & \cdot & \cdot & \cdot       & \id          & \cdot \\
		E     & -N    & \cdot & \cdot & \cdot       & \cdot        & \cdot \\
		\cdot & \cdot & \cdot & \cdot & \cdot       & \cdot        & \cdot \\
		\cdot & \cdot & \cdot & \cdot & \cdot       & \cdot        & \cdot \\
		\cdot & \cdot & \cdot & \cdot & E^\newtonupdate & \cdot        & E^\newtonupdate \\
		\cdot & \cdot & \cdot & \cdot & \cdot       & \cdot        & -(N^\newtonupdate)\adjoint \\
		\cdot & \cdot & \cdot & \cdot & E^\newtonupdate & -N^\newtonupdate & \cdot
	\end{pmatrix}
	\!
	\begin{pmatrix}
		\newtonupdate \\
		G \\
		u^\newtonupdate \\
		p^\newtonupdate \\
		\shapeprojgrad^\newtonupdate \\
		F^\newtonupdate \\
		\Pi^\newtonupdate
	\end{pmatrix}
	+
	\begin{pmatrix}
		\cdot \\
		\cdot \\
		\frac{\partial}{\partial u}             \LL( \newtonupdate, u^\newtonupdate, p^\newtonupdate) \\
		\frac{\partial}{\partial p}             \LL( \newtonupdate, u^\newtonupdate, p^\newtonupdate) \\
		\frac{\partial}{\partial \newtonupdate} \LL( \newtonupdate, u^\newtonupdate, p^\newtonupdate) \\
		\cdot \\
		\cdot
	\end{pmatrix}
	=
	0
	.
\end{equation}
As before, $N$ and $E$ denote the normal force operator and the elasticity operator on $\Omega$.

The system \eqref{eq:massive_nonlinear_system_2} for $W$ and the further, auxiliary unknowns corresponds to the nonlinear system $\Phi(x + \delta x) = 0$ for the step $\delta x$.
For convenience,
we recall the meaning
of the seven equations in \eqref{eq:massive_nonlinear_system_2}.
The first equation requires $F^\newtonupdate = 0$, i.e., the stationarity of the updated domain $\Omega_\newtonupdate$.
The second equation is the requirement that the displacement $\newtonupdate$ is induced by the (normal) force $G$.
The third and fourth equation are the adjoint and state equation on $\Omega_\newtonupdate$.
Finally, the last three equations
are the pull-back of the system \eqref{eq:continuous_saddle_point_system} on $\Omega_\newtonupdate$
to $\Omega$.

We can now describe a step of our Newton-like procedure for the solution of the nonlinear system \eqref{eq:massive_nonlinear_system_2}.
Suppose that $\Omega$ is the current domain and consider an iterate of the form $(0, 0, u, p, \shapeprojgrad, F, \Pi)$ with the state, the adjoint state, and the solution of \eqref{eq:continuous_saddle_point_system} on $\Omega$.
Notice that for this iterate, the residual of \eqref{eq:massive_nonlinear_system_2} is $(F,0,0,0,0,0,0)$.
Next we linearize the system \eqref{eq:massive_nonlinear_system_2} about this current iterate w.r.t.\ all seven variables.
We refrain from stating the lengthy formula for the linear system which results.
In practice, we generate this linear system governing the Newton step using the algorithmic differentiation capabilities of \fenics\ (\cite{LoggMardalWells2012:1}).
From the solution of that linear system we only extract the Newton update for the perturbation field.
We refer to it as $\newtonupdate$ since its current value is zero.
We then apply $\newtonupdate$ to the current domain $\Omega$ to obtain the new domain $(\id + \newtonupdate)(\Omega)$.
The six remaining variables are updated in a different fashion.
Rather than using the solution from the Newton step, we solve again the state and adjoint state equations on the new domain, as well as the system \eqref{eq:continuous_saddle_point_system} returning the projected shape gradient.
This procedure can be understood as a Newton-like method with nonlinear updates for some of the variables.
It ensures that the new iterate is of the same form as above.
Moreover, it allows us access to the projected shape gradient and its norm in every iteration so that we can use $\norm{\shapeprojgrad}_E \le \varepsilon_{\textup{tol}}$ as a stopping criterion as we did for the restricted gradient method.

Numerically, we have observed some instabilities if the current iterate $\Omega$ is far from being stationary.
Moreover we wish to establish a step size control in order to monitor the Armijo condition \eqref{eq:armijo_backtracking} and the mesh quality condition \eqref{eq:mesh_quality}.
To this end we added a regularization term $- G/\alpha$ to the first equation of \eqref{eq:massive_nonlinear_system_2}, i.e., we obtain
\begin{equation}
	\label{eq:massive_nonlinear_system_3}
	\begin{pmatrix}
		\cdot & -\alpha^{-1} \id       & \cdot & \cdot & \cdot           & \id              & \cdot \\
		E     & -N                     & \cdot & \cdot & \cdot           & \cdot            & \cdot \\
		\cdot & \cdot                  & \cdot & \cdot & \cdot           & \cdot            & \cdot \\
		\cdot & \cdot                  & \cdot & \cdot & \cdot           & \cdot            & \cdot \\
		\cdot & \cdot                  & \cdot & \cdot & E_\newtonupdate & \cdot            & E_\newtonupdate \\
		\cdot & \cdot                  & \cdot & \cdot & \cdot           & \cdot            & -N_\newtonupdate\adjoint \\
		\cdot & \cdot                  & \cdot & \cdot & E_\newtonupdate & -N_\newtonupdate & \cdot
	\end{pmatrix}
	\begin{pmatrix}
		\newtonupdate \\
		G \\
		u \\
		p \\
		\shapeprojgrad \\
		F \\
		\Pi
	\end{pmatrix}
	+
	\begin{pmatrix}
		\cdot \\
		\cdot \\
		\frac{\partial}{\partial u}             \LL( \newtonupdate, u, p ) \\
		\frac{\partial}{\partial p}             \LL( \newtonupdate, u, p ) \\
		\frac{\partial}{\partial \newtonupdate} \LL( \newtonupdate, u, p ) \\
		\cdot \\
		\cdot
	\end{pmatrix}
	=
	0
	.
\end{equation}
Thus, the update resulting from the solution of the Newton system satisfies
$- \alpha^{-1} \delta G + \delta F = - F$.
Heuristically,
this leads to $\delta G \approx \alpha \, F$ for small $\alpha$.
Consequently,
the transformation field 
which is applied to the current domain $\Omega$ satisfies
$\newtonupdate = E^{-1} \, N \,\delta G \approx \alpha \, \shapeprojgrad$
is essentially a scaled (restricted) gradient direction for small $\alpha$.
Therefore, similar as in a Levenberg--Marquardt method, we will refer to $\alpha$ as the damping parameter and it serves the same purpose as the step length parameter in \cref{alg:restricted_gradient_descent}.

A discrete variant of our Newton-like method is readily derived and given as \cref{alg:restricted_newton}.
In order to determine an appropriate damping parameter we consider analogues of the Armijo condition \eqref{eq:armijo_backtracking} and the mesh quality criterion \eqref{eq:mesh_quality}.
For the sake of clarity we re-state them with the relevant quantities for the Newton-like method.
In particular, we use the step length $\alpha = 1$ therein,
since the scaling of the step is already realized by the damping in \eqref{eq:massive_nonlinear_system_3}.
The Armijo condition becomes
\begin{equation}
	\label{eq:armijo_backtracking_newton}
	J_h\bigh(){ (\id + \newtonupdatedisc)(\Omega_h) }
	\le
	J_h( \Omega_h )
	+
	\sigma \, J_h'(\Omega_h; \newtonupdatedisc).
\end{equation}
with some parameter $\sigma \in (0,0.5)$.
The mesh quality criterion holds if
\begin{equation}
	\label{eq:mesh_quality_newton}
	\frac12 \le \det(\id + D\newtonupdatedisc) \le 2
	,
	\qquad
	\norm{ D\newtonupdatedisc }_F \le 0.3
	.
\end{equation}
is satisfied in every cell.
In addition we verify that $\newtonupdatedisc$ yields a descent direction.
If any of the above conditions fails, we decrease the damping parameter $\alpha$.

\begin{algorithm2e}[htp]
	\SetAlgoLined
	\KwData{Initial domain $\Omega_h$ \\
		Initial damping parameter $\alpha$,
		convergence tolerance $\varepsilon_{\textup{tol}}$,\\
		line search parameters $\beta \in (0,1)$, $\sigma \in (0,0.5)$
	}
	\KwResult{Improved domain $\Omega_h$ on which \eqref{eq:discrete_stationary_point_2} holds up to $\varepsilon_{\textup{tol}}$}
	\For{$i \leftarrow 1$ \KwTo $\infty$}{
		Solve the discrete state equation \eqref{eq:weak_form_example_PDE_discrete} for $u_h$\;
		Solve the discrete adjoint equation \eqref{eq:adjoint_problem_discrete} for $p_h$\;
		Solve \eqref{eq:discrete_shape_gradient_2} for $\shapeprojgraddisc$ with shape derivative $J_h'(\Omega_h;\cdot)$ from \eqref{eq:shape_derivative_discrete}\;
		\If{${\dual{E_h\shapeprojgraddisc}{\shapeprojgraddisc}} \le \varepsilon_{\textup{tol}}^2$}{
				STOP, the current iterate $\Omega_h$ is almost stationary for \eqref{eq:discrete_stationary_point_2}\;
			}
			Increase damping parameter $\alpha \leftarrow \alpha / \beta$\;
			Solve the Newton system associated with \eqref{eq:massive_nonlinear_system_3} with damping parameter $\alpha$ and extract the first component as $\newtonupdatedisc$\;
		\While{$J_h'(\Omega_h; \newtonupdatedisc) \ge 0$ holds, or \eqref{eq:armijo_backtracking_newton} or \eqref{eq:mesh_quality_newton} is violated}{
				Decrease damping parameter $\alpha \leftarrow \beta \, \alpha$\;
				Solve the Newton system associated with \eqref{eq:massive_nonlinear_system_3} with damping parameter $\alpha$ and extract the first component as $\newtonupdatedisc$\;
			}
			Transform the domain according to $\Omega_h \leftarrow (\id + \newtonupdatedisc)(\Omega_h)$\;
		}
	\caption{Restricted Newton method for \eqref{eq:discrete_stationary_point_2}.}
	\label{alg:restricted_newton}
\end{algorithm2e}

\section{Numerical Results: Newton-Like Method}
\label{sec:numerical_results_newton}

\begin{figure}
	\begin{subfigure}[b]{\linewidth}
		\centering
		\includegraphics[scale=0.15]{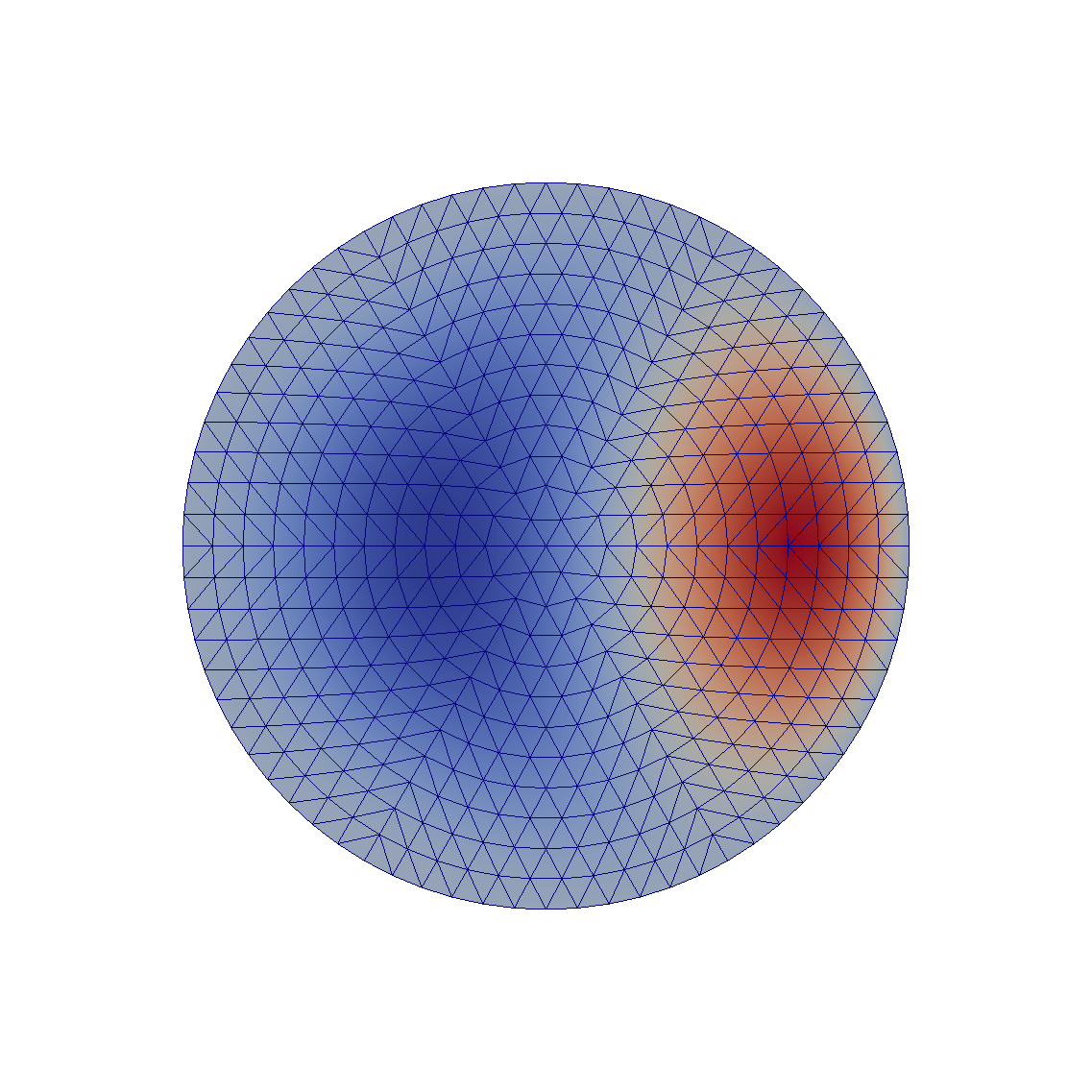}
		\hfill
		\includegraphics[scale=0.15]{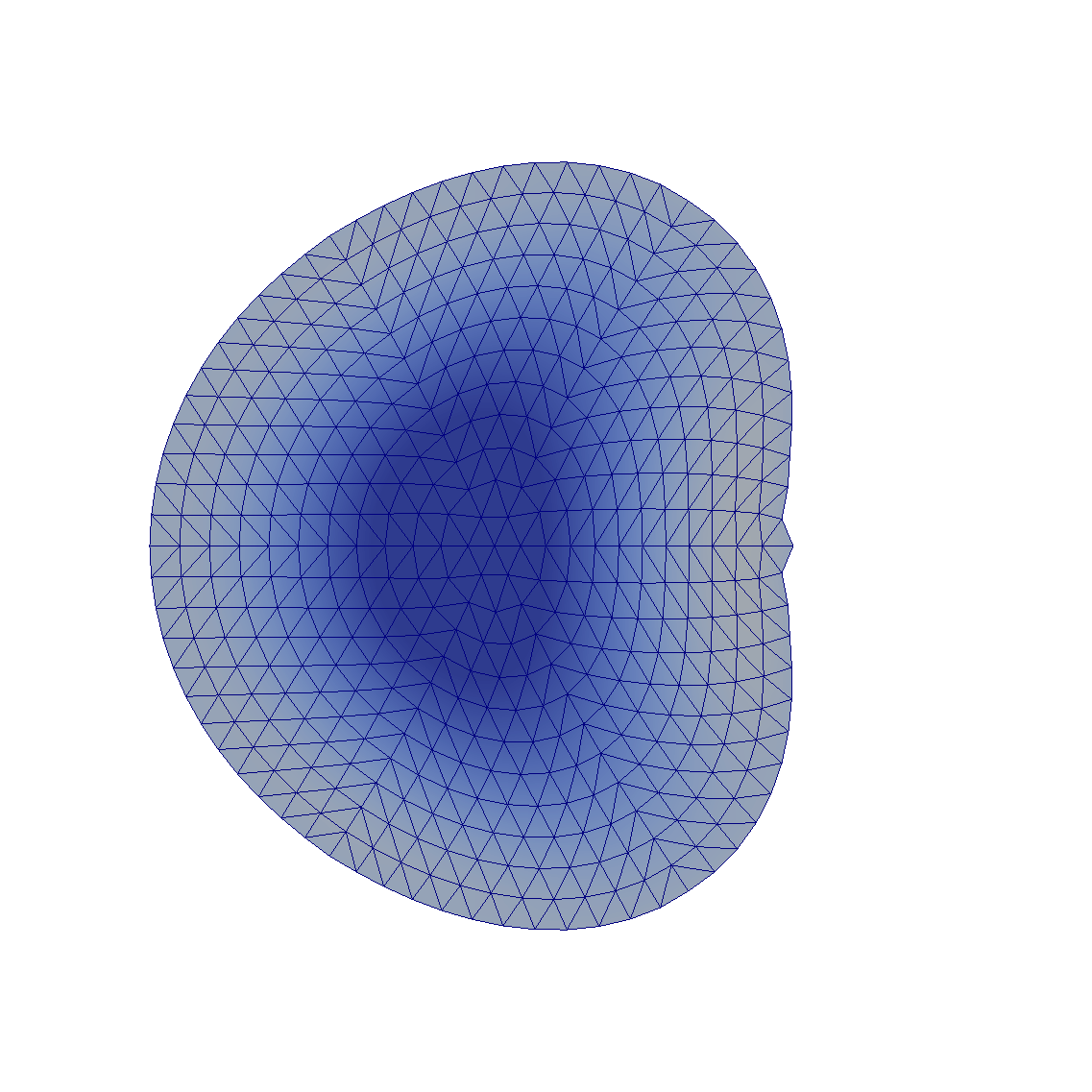}
	\end{subfigure}%

	\begin{subfigure}[b]{\linewidth}
		\centering
		\includegraphics[scale=0.15]{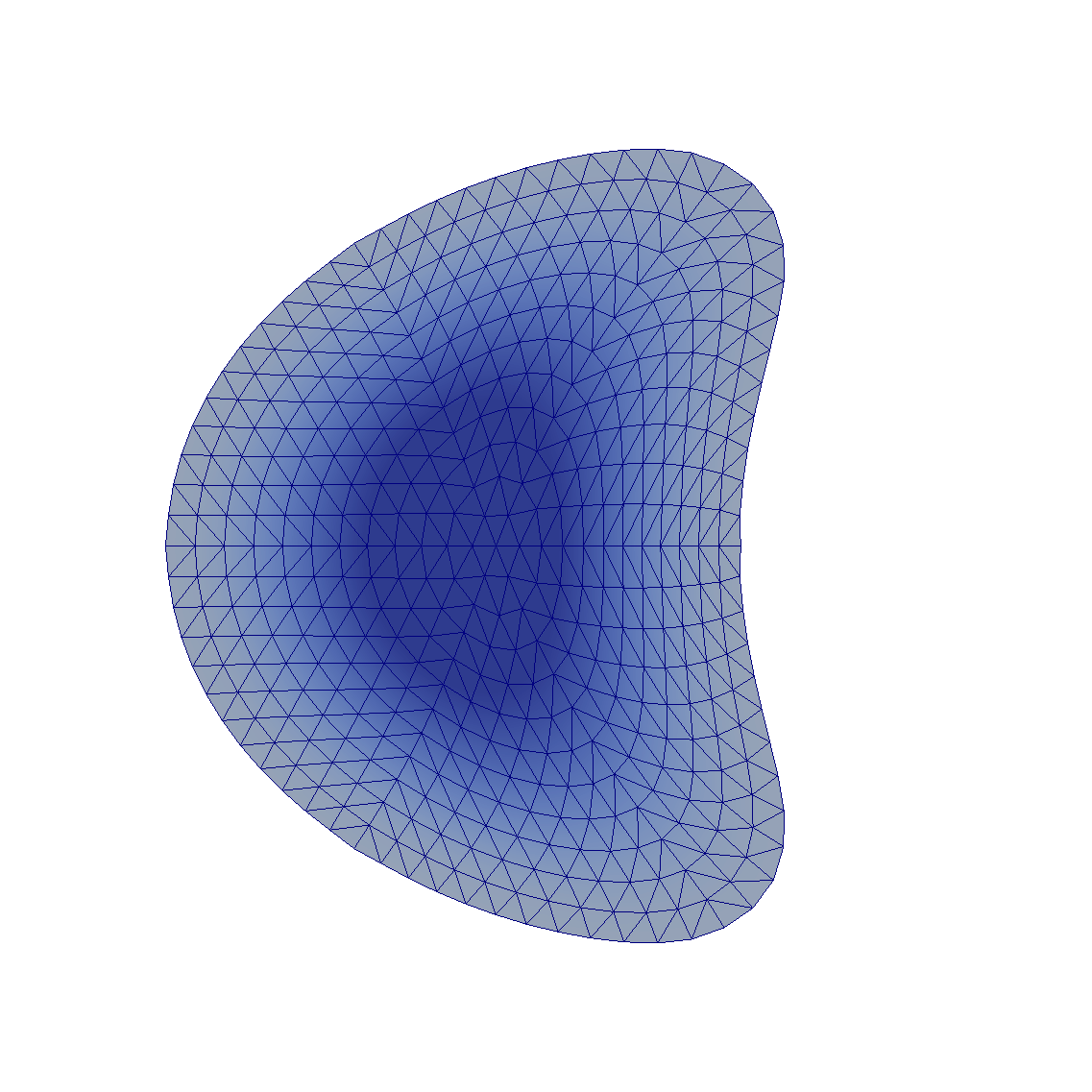}
		\hfill
		\includegraphics[scale=0.15]{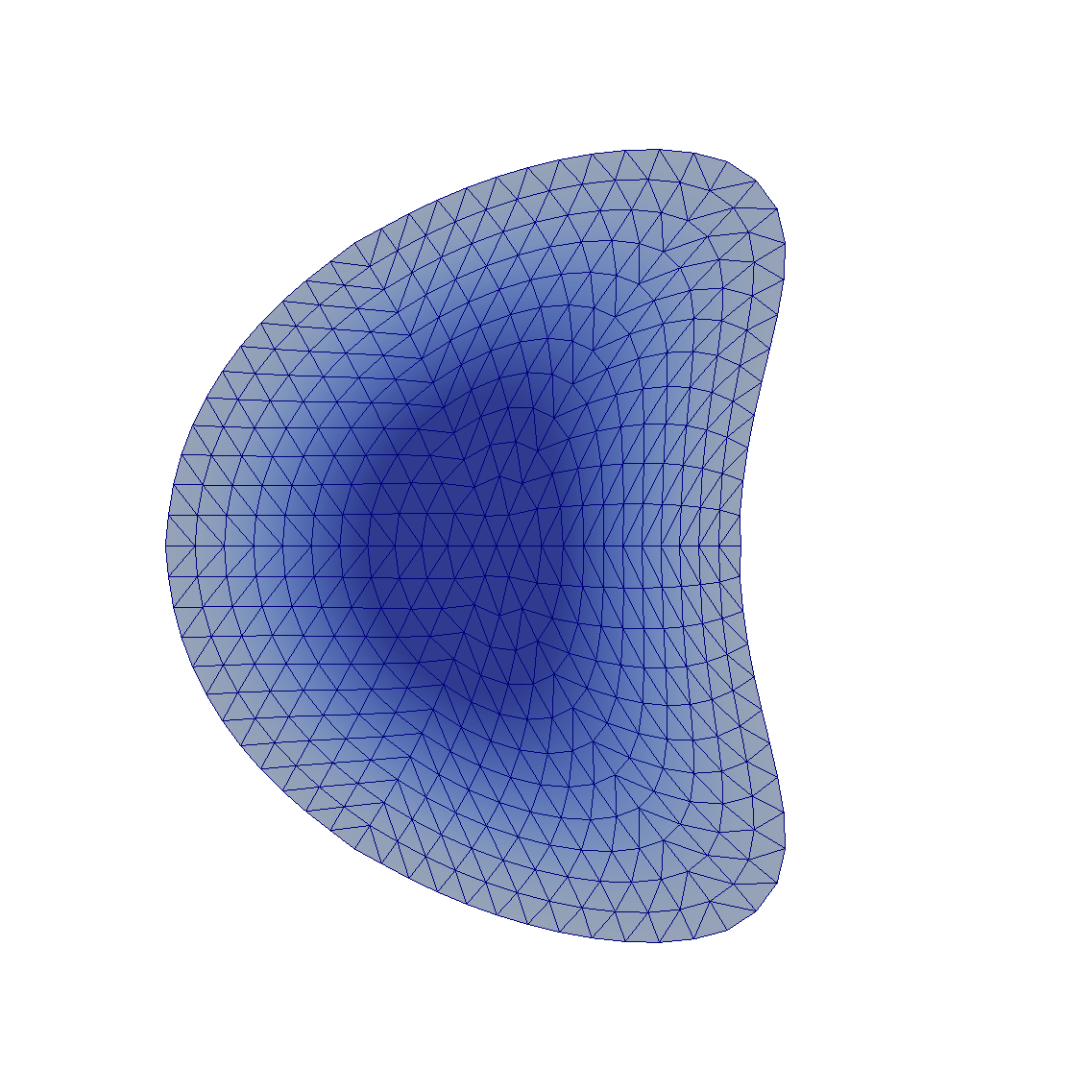}
	\end{subfigure}
	\caption{Intermediate shapes $\Omega_h$ obtained with the restricted Newton method at iterations~0, 4, 9, 12.}
	\label{fig:shapes_newton}
\end{figure}

This section is devoted to numerical results obtained by solving the same problem as in
\cref{sec:numerical_results_gradient} using the Newton-like method as described in 
\cref{sec:restricted_shape_hessian}. 
As mentioned in \cref{sec:numerical_results_gradient}, our implementation is freely available, see \cite{EtlingHerzogLoayzaWachsmuth2018:2}.
For this approach the stopping criterion 
\begin{equation}
	\label{eq:Newton_stopping_criterion}
	\norm{\shapeprojgraddisc}_{E_h}
	\le 
	\varepsilon_{\textup{tol}} = 10^{-9}
\end{equation}
was satisfied after 12~iterations and 6~seconds
on the previously used mesh with 469~vertices and 864~elements.  
In this case we used the line search parameters $\beta = 0.1$ and $\sigma = 0.1$ and an initial value of $\alpha = 10^{-2}$.
Young's modulus and the Poisson ratio are $E_0 = 1.0$ and $\nu = 0.4$ as before.
Some of the intermediate shapes are shown in \cref{fig:shapes_newton}.
As was already mentioned, we have the linear system in each Newton step
assembled using the algorithmic differentiation capabilities of \fenics\ and
solved in the same way as we did for the gradient method.
In this scenario the geometry condition \eqref{eq:mesh_quality_newton} led to a
decrease in the damping parameter $\alpha$ four times total (in iterations~3,
4, 5, and 7), while the Armijo condition \eqref{eq:armijo_backtracking_newton}
never necessitated a decrease in $\alpha$. We conjecture that the geometry
condition \eqref{eq:mesh_quality_newton} is triggered here more often compared
to the gradient methods since the Newton-like method tends to produce steps
with larger norm $\norm{W_h}_{E_h}$. However the perturbation of identity
transformation imposes a limit on the step size, which leads to a reduction of
$\alpha$. In the final iteration, the value of $\alpha = 10^5$ is reached.

	The experiments up to here were obtained on a coarse mesh with 469~vertices and 864~elements. 
	We also studied the dependence of iteration numbers and CPU time on the mesh level, both for the restricted shape gradient and restricted Newton methods.
	Finer mesh levels are obtained by uniform refinement.
	\Cref{tab:mesh_level} shows the size of the mesh in terms of the number of cells and vertices together with the number 
	of iterations required for the convergence of both algorithms, and the time of execution.

\begin{table}
	\begin{center}
		\begin{tabular}{|r|r|r|r|r|r|r|}
			\hline
			\multicolumn{2}{|c|}{\textbf{Mesh Level}}                            & \multicolumn{2}{c|}{\textbf{Restricted Gradient}}                      & \multicolumn{2}{c|}{\textbf{Restricted Newton}}                        \\ \hline
			\multicolumn{1}{|c|}{\textbf{Vertices}} & \multicolumn{1}{c|}{\textbf{Cells}} & \multicolumn{1}{c|}{\textbf{Iter}} & \multicolumn{1}{c|}{\textbf{Time {[}s{]}}} & \multicolumn{1}{c|}{\textbf{Iter}} & \multicolumn{1}{c|}{\textbf{Time {[}s{]}}} \\ \hline
			127                         & 216                           & 527                       & 10                                & 9                         & 3                                 \\ \hline
			469                         & 864                           & 864                       & 38                                & 11                        & 7                                 \\ \hline
			1801                        & 3456                          & 1481                      & 244                               & 13                        & 48                                \\ \hline
			7057                        & 13824                         & 2353                      & 1733                              & 14                        & 319                               \\ \hline
		\end{tabular}
		\caption{Number of iterations and time of execution for the 2D~example at different mesh levels to reach the tolerance~\eqref{eq:restricted_gradient_stopping_criterion} for the restricted gradient method. Moreover, for the restricted Newton method we used the stopping criterion~\eqref{eq:Newton_stopping_criterion} with a tolerance of $10^{-8}$ and an initial damping parameter $\alpha_0 = 10^{7}$.}
		\label{tab:mesh_level}
	\end{center}
\end{table}

\begin{figure}
	\begin{subfigure}[b]{\linewidth}
		\centering
		\includegraphics[scale=0.16]{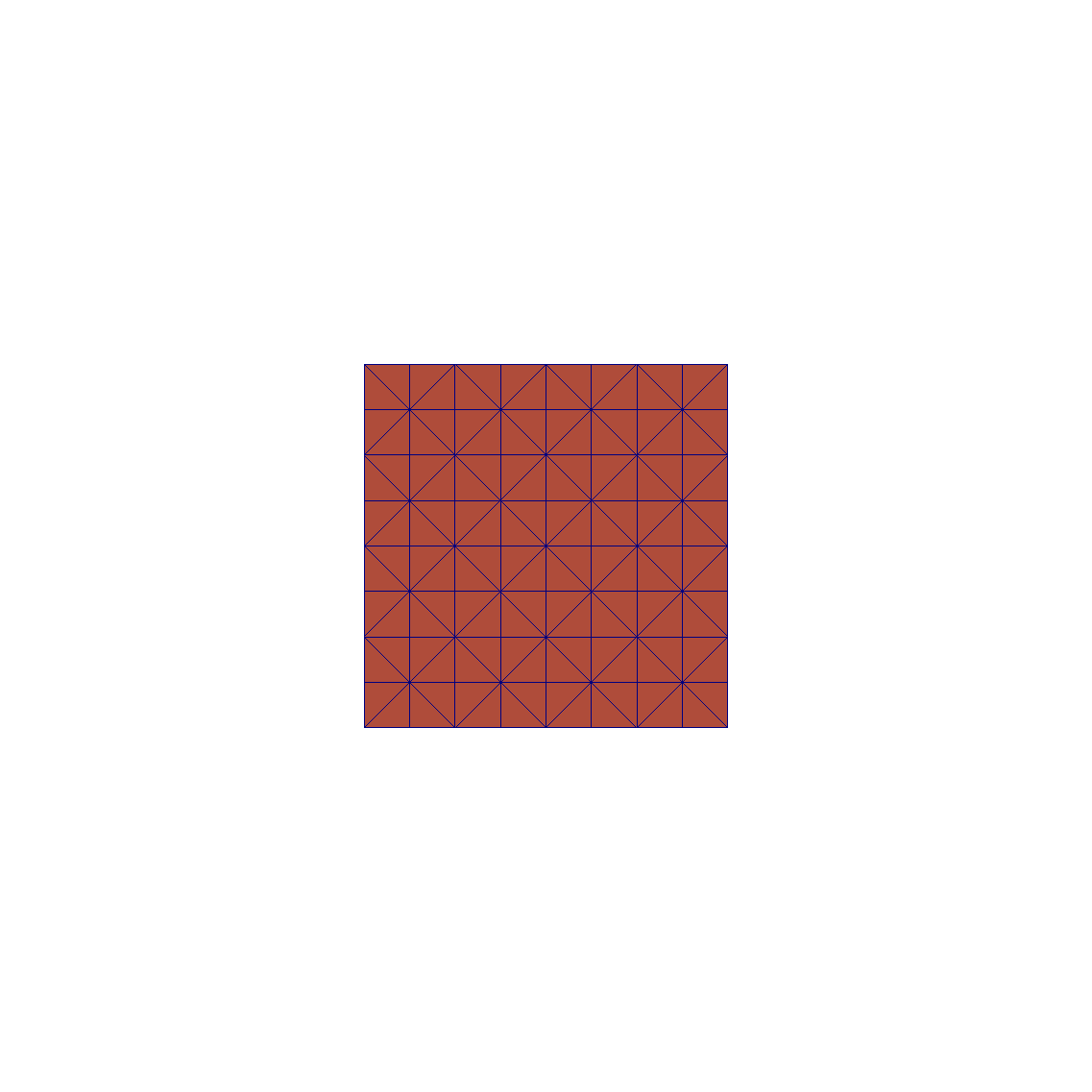}
		\hfill
		\includegraphics[scale=0.16]{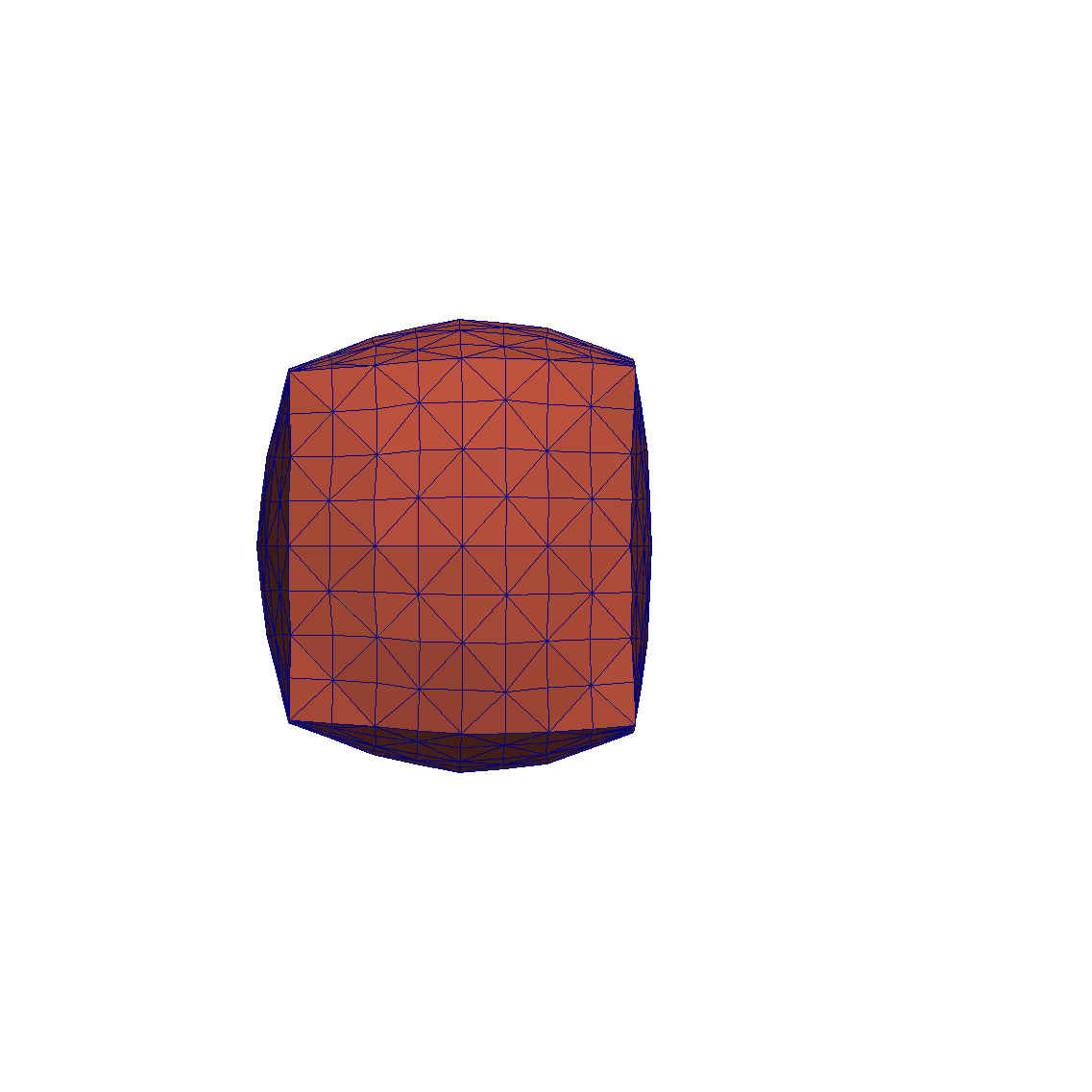}
	\end{subfigure}%
	\\[-.5cm]
	\begin{subfigure}[b]{\linewidth}
		\centering
		\includegraphics[scale=0.16]{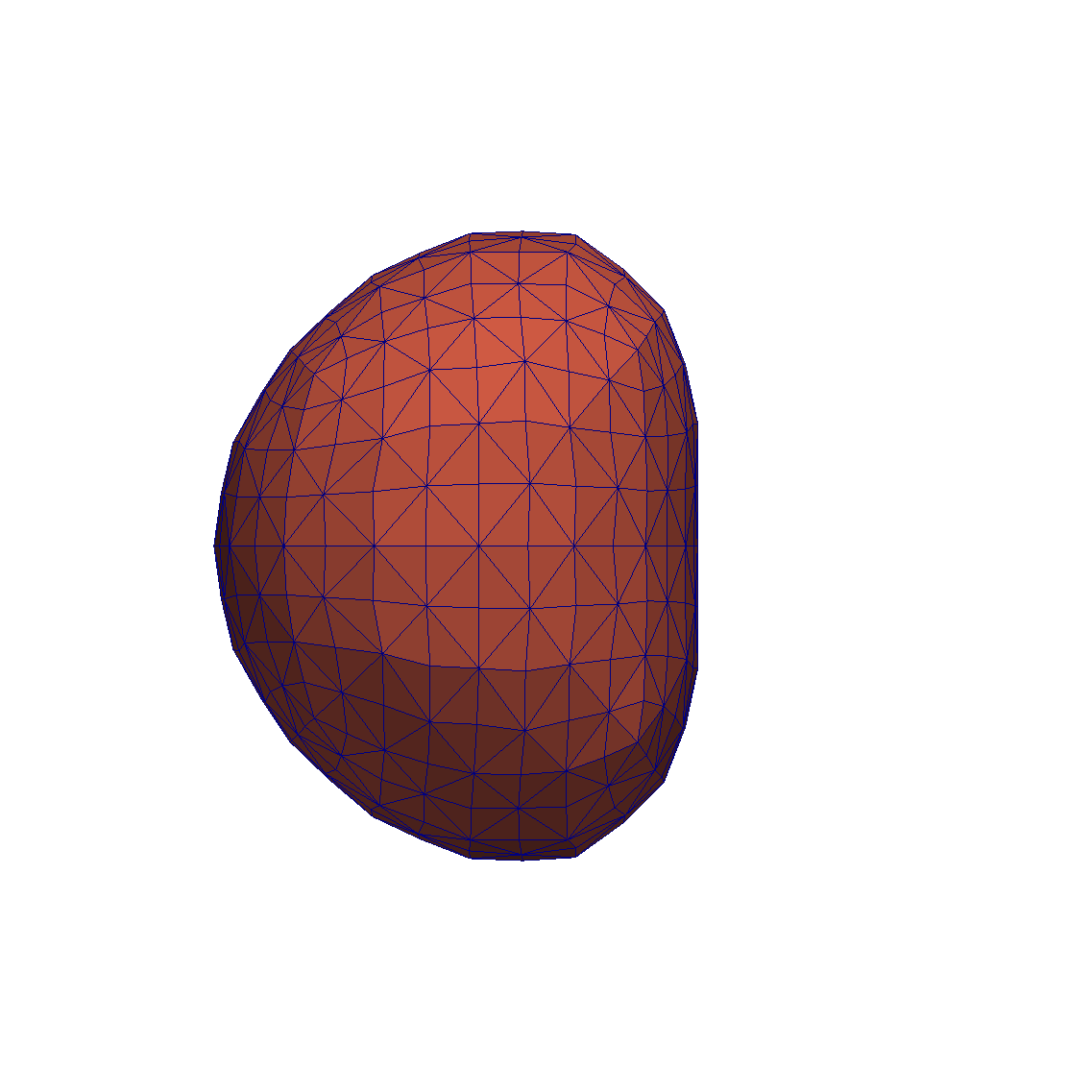}
		\hfill
		\includegraphics[scale=0.16]{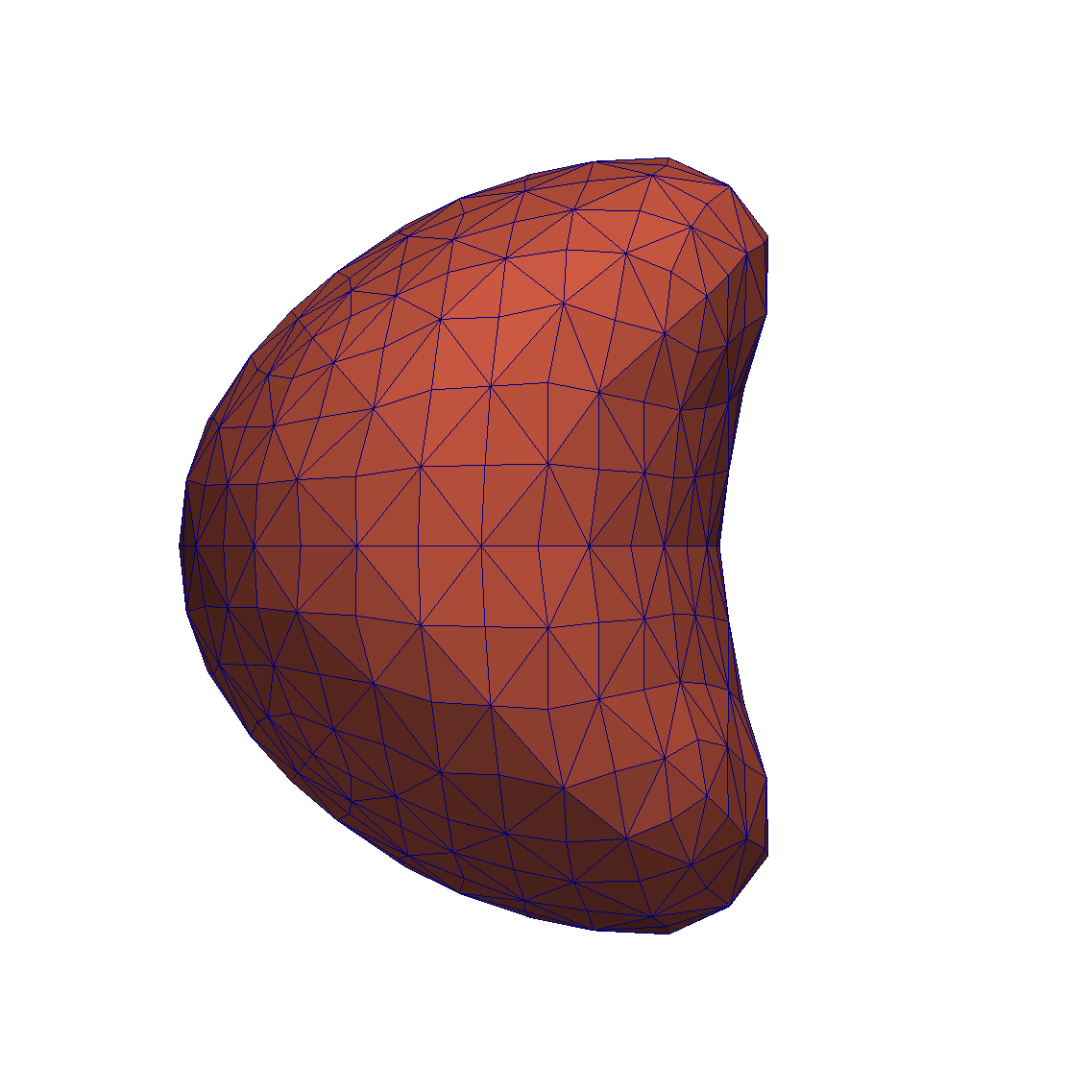}
	\end{subfigure}%
	\caption{Intermediate shapes $\Omega_h$ obtained with the restricted Newton method at iterations~0, 7, 12, 21 for the 3D example.}
	\label{fig:shapes_newton3D}
\end{figure}

Finally, we provide the results of the Newton-like method for an example in 3D.
This time, the right hand side of the state equation in \eqref{eq:shape_optimization_problem_continuous} is given by $f(x,y,z) = 2.5 \, (x + 0.4 - y^2)^2 + x^2 + y^2 + z^2 - 1$.
Moreover, the initial shape is a cube with 729~vertices and 3072~elements.
All other data remains the same as above.
\Cref{fig:shapes_newton3D} shows the initial, some intermediate and the final shapes obtained using~\cref{alg:restricted_newton} after 21~iterations and 286~seconds using the same tolerance as in the 2D case. 
Concerning the Armijo and geometry conditions \eqref{eq:armijo_backtracking_newton} and \eqref{eq:mesh_quality_newton} and the value of the damping parameter $\alpha$, we observe a very similar behavior as in the 2D case. The final value is $\alpha = 10^6$.

\section{Conclusions}
\label{sec:conclusions}

In this paper we introduce the concept of restricted mesh deformations for the computational solution of shape optimization problems involving PDEs.
In a nutshell, we only admit perturbations fields which are induced by normal boundary forces.
We argue that the stationarity condition \eqref{eq:discrete_stationary_point} which does not impose any restriction on the mesh deformations leads to degenerate meshes and premature stopping.
By contrast, we were able to solve the corresponding restricted stationarity condition \eqref{eq:discrete_stationary_point_2} to high accuracy even with a gradient method.
We also propose a Newton-like method based on restricted mesh deformations which exhibits fast convergence.

It is not clear whether \eqref{eq:discrete_stationary_point_2} are the optimality conditions of a discrete optimization problem in Euclidean space.
We conjecture that \eqref{eq:discrete_stationary_point_2} are the optimality conditions for a problem defined on a discrete shape manifold, whose tangent space is represented by restricted mesh deformations.

\section*{Acknowledgments}

The third author acknowledges support by the German Academic Exchange Service (DAAD) within the Doctoral Programm~2017/18.

\printbibliography
\end{document}